\RequirePackage{ifpdf}
\ifpdf 
\documentclass[pdftex]{sigma}
\else
\documentclass{sigma}
\fi

\numberwithin{equation}{section}

\begin{document}

\allowdisplaybreaks

\renewcommand{\thefootnote}{$\star$}

\renewcommand{\PaperNumber}{069}

\FirstPageHeading

\ShortArticleName{Homogeneous Poisson Structures on Loop Spaces of
Symmetric Spaces}

\ArticleName{Homogeneous Poisson Structures\\ on Loop Spaces of
Symmetric Spaces\footnote{This paper is a
contribution to the Special Issue on Kac--Moody Algebras and Applications. The
full collection is available at
\href{http://www.emis.de/journals/SIGMA/Kac-Moody_algebras.html}{http://www.emis.de/journals/SIGMA/Kac-Moody{\_}algebras.html}}}

\AuthorNameForHeading{D.~Pickrell}

\Author{Doug PICKRELL}

\Address{Department of Mathematics, University of Arizona, Tucson,
AZ, 85721, USA}

\Email{\href{mailto:pickrell@math.arizona.edu}{pickrell@math.arizona.edu}}

\ArticleDates{Received June 14, 2008, in f\/inal form September
27, 2008; Published online October 07, 2008}

\Abstract{This paper is a sequel to [Caine A., Pickrell D., {\it Int. Math. Res. Not.},  to appear, \href{http://arxiv.org/abs/0710.4484}{arXiv:0710.4484}], where we
studied the Hamiltonian systems which arise from the Evens--Lu
construction of homogeneous Poisson structures on both compact and
noncompact type symmetric spaces. In this paper we consider loop
space analogues. Many of the results extend in a
relatively routine way to the loop space setting, but new issues
emerge. The main point of this paper is to spell out the meaning
of the results, especially in the $SU(2)$ case. Applications
include integral formulas and factorizations for Toeplitz
determinants.}


\Classification{22E67; 53D17; 53D20}

\vspace{-1mm}

\renewcommand{\thefootnote}{\arabic{footnote}}
\setcounter{footnote}{0}

\section{Introduction}\label{Introduction}

The f\/irst purpose of this paper is to generalize the framework in
\cite{CP} to loop spaces. This generalization is straightforward,
using the fundamental insight of Kac and Moody that f\/inite
dimensional complex semisimple Lie algebras and (centrally
extended) loop algebras f\/it into the common framework of Kac--Moody
Lie algebras.

Suppose that $\dot {X}$ is a simply connected compact symmetric
space with a f\/ixed basepoint. From this, as we will more fully
explain in Sections \ref{loopgroups} and \ref{typeIcase}, we
obtain a diagram of groups
\begin{gather}\label{groupdiagram}\begin{matrix} &&G=\widehat {L}\dot {G}\\[-3pt]
&\nearrow&&\nwarrow\\[-3pt]
G_0=\widehat {L}\dot {G_0}&&&&U=\widehat {L}\dot {U}\\[-3pt]
&\nwarrow&&\nearrow\\[-3pt]
&&K=\widehat {L}\dot {K}\end{matrix}  \end{gather} where $\dot
{U}$ is the universal covering of the identity component of the
isometry group of $\dot {X}$, $\dot {X}\simeq \dot {U}/\dot {K}$,
$\dot {G}$ is the complexif\/ication of $\dot {U}$, $\dot {X}_0=\dot
{G}_0/\dot {K}$ is the noncompact type symmetric space dual to
$\dot {X}$, $L\dot {G}$ denotes the loop group of $\dot {G}$,
$\widehat {L}\dot {G}$ denotes a Kac--Moody extension, and so on.
This diagram is a prolongation of diagram (0.1) in \cite{CP}
(which is embedded in (\ref{groupdiagram}) by considering constant
loops).

We also obtain a diagram of equivariant totally geodesic (Cartan)
embeddings of symmetric spaces:
\begin{gather}\label{Cartandiagram}\begin{matrix} L(\dot {U}/\dot {K})&\stackrel{\phi}{\rightarrow}&\widehat {L} \dot {U}\\[-1pt]
\downarrow&&\downarrow\\[-1pt]
\tilde{L}\dot {G}/\tilde{L}\dot {G}_0&\stackrel{\phi}{\rightarrow}&\widehat {L} \dot {G}&
\stackrel{\psi}{\leftarrow}& \tilde{L}\dot {G}/\tilde{L}\dot {U}\\[-1pt]
&&\uparrow&&\uparrow\\[-1pt]
&&\widehat {L} \dot{G}_0&\stackrel{\psi}{\leftarrow}&L(\dot
{G}_0/\dot {K})\end{matrix}\end{gather} This is a prolongation
of diagram (0.2) in \cite{CP}.

Let $\Theta$ denote the involution corresponding to the pair
$(\dot {U},\dot {K})$. We consider one additional ingredient: a
triangular decomposition
\begin{gather}\label{triangularstructure}\dot {\mathfrak g}=\dot {\mathfrak n}^- \oplus \dot {\mathfrak h} \oplus \dot
{\mathfrak n}^+ \end{gather}which is $\Theta$-stable and for
which $\dot {\mathfrak t}_0=\dot {\mathfrak h}\cap \dot {\mathfrak
k}$ is maximal abelian in $\dot {\mathfrak k}$. There is a
corresponding Kac--Moody triangular decomposition
\[\widehat {L}_{\rm pol}\dot {\mathfrak g}=\Big(\bigoplus_{n<0}\dot {\mathfrak g}z^n\oplus \dot {\mathfrak n}^{
-}\Big)\oplus\mathfrak h\oplus \Big(\dot {\mathfrak
n}^{+}\oplus\bigoplus_{n>0}\mathfrak g z^n\Big)\] extending
(\ref{triangularstructure}).

This data determines standard Poisson Lie group structures,
denoted $\pi_U$ and $\pi_{G_0}$, for the groups $U=\widehat {L}
\dot {U}$ and $G_0=\widehat {L} \dot {G}_0$, respectively.  By a
general construction of Evens and Lu~\cite{EL}, the symmetric
spaces $X= L \dot {X}$ and $X_0= L\dot {X}_0$ acquire Poisson
structures $\Pi_X$ and $\Pi_{X_0}$, respectively, which are
homogeneous for the respective actions of the Poisson Lie groups
$(U,\pi_U)$ and $(G_0,\pi_{G_0})$. These spaces are inf\/inite
dimensional, and there are many subtleties associated with Poisson
structures in inf\/inite dimensions (see \cite{OR}). Consequently in
this paper we will always display explicit decompositions and
formulas, and we will avoid any appeal to general theory (for
``symplectic foliations'', for example).

The plan of this paper is the following. In Section
\ref{loopgroups} we introduce notation and recall some well-known
facts concerning loop algebras and groups.

In Section \ref{typeIcase} we consider the case when $\dot {X}$ is
an irreducible type I space. All of the results of Sections 1--4 of
\cite{CP} generalize in a relatively straightforward way to the
loop context roughly outlined above. The basic result is that
$\Pi_{X_0}$ has just one type of symplectic leaf, this leaf is
Hamiltonian with respect to the natural action of $T_0$, there are
relatively explicit formulas for this Hamiltonian system, and in a
natural way, this system is isomorphic to the generic Hamiltonian
system for $\Pi_X$. Although this system is inf\/inite dimensional,
a heuristic application of the Duistermaat--Heckman exact
stationary phase theorem to this system suggests some remarkable
integral formulas. This is discussed in Section 7 of
\cite{pickrell1}. These formulas remain conjectural.

In Section \ref{S^2case} I have attempted to do some calculations
in the $\dot {X}=S^2$ case. The formulas are complicated; I
included them to give the reader a concrete feeling for the
subject.

In Section \ref{typeIIcase}, and Appendix~\ref{AppendixA}, we consider the group
case. Again, the results of Sec\-tions 1--4 of \cite{CP} generalize
in a straightforward way. However signif\/icant issues emerge when
we try to generalize the results of Section 5 of \cite{CP}. In the
f\/inite dimensional context of $\dot {X}=\dot {K}$, the (negative
of the) standard Poisson Lie group structure $(\dot {K}, \pi_{\dot
{K}})$ is isomorphic to $(\dot {X}, \Pi_{\dot {X}})$, by left
translation by a representative for the longest Weyl group
element. In the loop context the Poisson Lie group and Evens--Lu
structures are fundamentally dif\/ferent: the symplectic leaves for
$\pi_K$ (essentially Bruhat cells) are f\/inite dimensional, whereas
the symplectic leaves for $\Pi_X$ (essentially Birkhof\/f strata)
are f\/inite codimensional. In f\/inite dimensions Lu has completely
factored the symplectic leaves. Lu's results, as formulated in
\cite{Lu} in terms of $\pi_K$, do generalize in a relatively
straightforward way to the loop context. Some details of this
generalization are worked out in Appendix A, where we have
extended this to the larger category of symmetrizable Kac--Moody
algebras.

The basic question is whether the Hamiltonian systems for $\Pi_X$,
in this inf\/inite dimensional context, are solvable (in a number of
senses). In Sections \ref{typeIIcase}--\ref{SU(2)caseII} we show
that there is a natural way to conjecturally reformulate and
extend Lu's results to suggest that the generic symplectic leaves
are integrable. However we have not succeeded in fully proving
this conjecture: while we can factor the momentum mapping, and the
Haar measure relevant for Theorem \ref{hankeltheorem} below, we
have not shown that the symplectic form $\Pi_X^{-1}$ factors.

In Section \ref{SU(2)case} we spell out the meaning of the results
in Section \ref{typeIIcase} when $\dot {K}=SU(2)$. One consequence
is the following integral formula.

\begin{theorem}\label{hankeltheorem} Given $x_j\in \mathbb C$, let
\begin{gather}\label{hank1} B\left(\sum_{j=1}^nx_jz^j\right)=\left(\begin{array}{ccccc} x_n & 0 & & \dots & 0\\
x_{n-1}&x_n&0&\dots& 0\\
\vdots & &\ddots &\ddots & \vdots\\
x_2&\dots&x_{n-1}&x_n&0\\
x_1&x_2&\dots&x_{n-1}&x_n\end{array} \right).
 \end{gather} Then
\begin{gather*}
\int\prod_{l=0}^{n-1}\det\left(1+B\left(\sum_{j=1}^{n-l}x_{l+j}z^j\right)\right)B\left(\sum_{
j=1}^{n-l}\big(x_{l+j}z^j\big)^{*}\right)^{-p_l}d\lambda (x_1,\dots ,x_n)\\
\qquad{}=\pi^n\frac 1{(p_1-1)}\frac 1{(2p_1+p_2-3)}\cdots \frac 1{(np_1+(n-1
)p_2+\cdots +p_n-(2n-1))}.
\end{gather*} In particular, if we write $B_n(x)$ for the
matrix \eqref{hank1}, for a general power series $x=\sum x_jz^j$,
then
\begin{gather}\label{hank2}\frac
1{\det(1+B_n(x)B_n(x)^{*})^p}d\lambda (x_1,\dots,x_n)\end{gather} is
a finite measure if and only if $p>2-1/n$.
\end{theorem}

This result is important because it determines the critical
exponents for the integrands in~(\ref{hank2}) exactly, whereas I
am not aware of any other way to even estimate these exponents in
a~useful way.  The relevance of this to the theory of conformally
invariant measures, where one must understand the limit as
$n\to\infty$, is described in \cite{pickrell2}.

In Section~\ref{SU(2)caseII} we consider the question of global
solvability of the symplectic leaves, in the $SU(2)$ case. A
consequence of the global factorization of the momentum mapping is
the following illustrative statement about block Toeplitz
operators.

\begin{theorem} Given complex numbers $\eta_j$, $\chi_j$, $\zeta_j$, let
$g:S^1 \to SU(2)$ be the product of $SU(2)$ loops
\begin{gather*}
a(\eta_0)\left(\begin{matrix} 1&\eta_0\\
-\bar{\eta}_0&1\end{matrix}
\right)\cdots a(\eta_n)\left(\begin{matrix} 1&\eta_nz^n\\
-\bar{\eta}_nz^{-n}&1\end{matrix} \right) \left(\begin{matrix}
e^{\sum\chi_jz^j}&0\\0&e^{-\sum\chi_jz^j} \end{matrix}\right),\\
a(\zeta_n)\left(\begin{matrix} 1&\zeta_nz^{-n}\\
-\bar{\zeta}_nz^n&1\end{matrix}
\right)\cdots a(\zeta_1)\left(\begin{matrix} 1&\zeta_1z^{
-1}\\
-\bar{\zeta}_1z&1\end{matrix} \right),
\end{gather*}
 where $a(\cdot)=(1+\vert
\cdot\vert^2)^{-1/2}$ and $\chi_{-j}=-\bar{\chi_j}$. Let $A(g)$
denote the Toeplitz operator defined by the symbol $g$. Then
\[\det(A(g)A(g)^*)=\prod_ja(\eta_j)^{2j}a(\zeta_j)^{2j}e^{-\vert j
\vert \vert \chi_j \vert^2}.\]
\end{theorem}

When $\eta$ and $\zeta$ vanish, this reduces to a well-known
formula with a long history (e.g.\ see Theorem 7.1 of \cite{W}).

In Section \ref{SU(2)caseII}, because the $SU(2)$ loop space is
inf\/inite dimensional, it is necessary to take a~limit as
$n\to\infty$, so that the above product of loops is to be
interpreted as an inf\/inite factorization of a generic $g\in
LSU(2)$. At a heuristic level, the invariant measures considered
in \cite{pickrell1} factor in these coordinates. The conjectural
integral formulas in Section 7 of \cite{pickrell1} (in the $SU(2)$
case) follow immediately from this product structure. However
changing coordinates in inf\/inite dimensions is nontrivial, and
probabilistic analysis is required to justify this claim.

\section{Loop groups}\label{loopgroups}

In this section we recall how (extended) loop algebras f\/it into
the framework of Kac--Moody Lie algebras. The relevant structure
theory for loop groups is developed in \cite{PS}, and for loop
algebras in Chapter 7 of \cite{K}.

Let $\dot {U}$ denote a simply connected compact Lie group. To
simplify the exposition, we will assume that $\dot {\mathfrak u}$
is a simple Lie algebra. Let $\dot {G}$ and $\dot {\mathfrak g}$
denote the complexif\/ications, and f\/ix \mbox{a~$\dot {\mathfrak
u}$-compatible} triangular decomposition
\begin{gather}\label{6.1}\dot {\mathfrak g}=\dot {\mathfrak n}^{-}\oplus\dot {\mathfrak h}\oplus\dot {\mathfrak n}^{
+}.\end{gather} We let $\langle\cdot ,\cdot\rangle$ denote the
unique $Ad(\dot {G})$-invariant symmetric bilinear form such that
(for the dual form) $\langle\theta ,\theta\rangle = 2$, where
$\theta$ denotes the highest root for $\dot {\mathfrak g}$, i.e.\
$\langle\cdot ,\cdot\rangle=\frac 1{\dot {g}}\kappa$, where
$\kappa$ denotes the Killing form, and $\dot {g}$ is the dual
Coxeter number.

Let $\widehat {L}\dot {\mathfrak g}$ denote the real analytic
completion of the untwisted af\/f\/ine Lie algebra corresponding to
$\dot {\mathfrak g}$, with derivation included (the degree of
smoothness of loops is essentially irrelevant for the purposes of
this discussion; any f\/ixed degree of Sobolev smoothness $s>1/2$
would work equally well). This is def\/ined in the following way. We
f\/irst consider the universal central extension of $L\dot
{\mathfrak g}=C^{\omega}(S^1,\dot {\mathfrak g})$,
\[0\to \mathbb C c\to\tilde {L}\dot {\mathfrak g}\to L\dot {\mathfrak g}\to 0.\]
As a vector space $\tilde {L}\dot {\mathfrak g}=L\dot {\mathfrak
g}\oplus \mathbb C c$. In these coordinates, the $\tilde {L}\dot
{\mathfrak g}$-bracket is given by
\begin{gather}\label{bracket}[X+\lambda c,Y+\lambda^{\prime} c]_{\tilde {L}\dot {\mathfrak g}}=[X,Y]_{L\dot {\mathfrak g}}+\frac
i{2\pi}\int_{S^1}\langle X\wedge dY\rangle c. \end{gather} Then
$\widehat {L}\dot {\mathfrak g}=\mathbb C d\propto\tilde {L}\dot
{\mathfrak g}$ (the semidirect sum), where the derivation $d$ acts
by $d(X+\lambda c)=\frac 1i\frac d{d\theta}X$, for $X\in L\dot
{\mathfrak g}$. The algebra generated by $\dot {\mathfrak
u}$-valued loops induces a central extension
\[0\to i\mathbb R c\to\tilde {L}\dot {\mathfrak u}\to L\dot {\mathfrak u}\to 0 \]
and a real form $\widehat {L}\dot {\mathfrak u}=i\mathbb R
d\propto\tilde {L}\dot {\mathfrak u}$ for $ \widehat {L}\dot
{\mathfrak g}$. We identify $\dot {\mathfrak g}$ with the constant
loops in $L\dot {\mathfrak g}$. Because the extension is trivial
over $\dot {\mathfrak g}$, there are embeddings of Lie algebras
\[\dot {\mathfrak g}\to\tilde {L}\dot {\mathfrak g}\to\widehat {L}\dot {\mathfrak g}. \]

The Lie algebra $\widehat {L}\dot {\mathfrak g}$ has a triangular
decomposition
\begin{gather}\label{looptriangulardecomposition}\widehat {L}\dot {\mathfrak g}
=\mathfrak n^{-}\oplus \mathfrak h \oplus \mathfrak n^{ +}, \end{gather} where
\begin{gather*}
 \mathfrak h=\dot {\mathfrak h}+\mathbb C c+\mathbb C d, \\
\mathfrak n^{+}=\left\{x=\sum_0^{\infty}x_nz^n\in H^0(D;\dot {\mathfrak g}):x(0)=x_0\in\dot {\mathfrak n}^{+}\right\}
\end{gather*}
and
\[
\mathfrak n^{-}=\left\{x=\sum_0^{\infty}x_nz^{-n}\in H^0(D^{*};\dot {\mathfrak g}):x(\infty )=x_0\in\dot {\mathfrak n}^{
-}\right\}.
\] This is compatible with the f\/inite dimensional triangular
decomposition (\ref{6.1}).  We let $N^{\pm}$ denote the prof\/inite
nilpotent groups corresponding to $\mathfrak n^{\pm}$, e.g.
\[N^-=H^0(D^*,\infty;\dot{G},\dot{N}^-).\]

There is a unique ${\rm Ad}$-invariant symmetric bilinear form on
$\widehat {L}\dot {\mathfrak g}$ which extends the normalized
Killing form on $\dot {\mathfrak g}$. It has the following
restriction to $\mathfrak h$:
\[\langle c_1d+c_2c+h,c_1'd+c_2'c+h'\rangle
=c_1c_2'+c_2c_1'+\langle h,h'\rangle .\] This form is
nondegenerate. The restriction of this form to $\widehat {L} \dot
{u}$ is also nondegenerate, although this restriction is of
Minkowski type, in contrast to the f\/inite dimensional situation.

The simple roots for $(\widehat {L}\dot {\mathfrak g},\mathfrak
h)$ are $\{\alpha_j: 0\le j\le rk\dot {\mathfrak g}\}$, where
\[\alpha_0=d^*-\theta ,\qquad\alpha_j=\dot{\alpha}_j,\qquad j>0, \]
$d^*(d)=1$, $d^*(c)=0$, $d^* (\dot {\mathfrak h})=0$, and the
$\dot{\alpha}_j$ denote the simple roots for the triangular
decomposition of $\dot {\mathfrak g}$ (with $\dot{\alpha}_ j$
vanishing on $c$ and $d$). The simple coroots of $\mathfrak h
\subset\widehat {L}\dot {\mathfrak g}$ are $\{h_j:0\le j\le rk\dot
{\mathfrak g}\}$, where
\[h_0=c-\dot {h}_{\theta},\qquad h_j=\dot {h}_j,\qquad j>0,\]
and the $\{\dot {h}_j\}$ are the simple coroots of $\dot
{\mathfrak g}$. For $i>0$, the root homomorphism $i_{\alpha_i}$ is
$i_{\dot{\alpha}_i}$ followed by the inclusion $\dot {\mathfrak
g}\subset\widehat {L}\dot {\mathfrak g}$.  For $i=0$
\begin{gather}\label{roothom}i_{\alpha_0}\left(\left(\begin{matrix} 0&0\\
1&0\end{matrix} \right)\right)=e_{\theta}z^{-1},\qquad
i_{\alpha_0}\left(\left(\begin{matrix}
0&1\\
0&0\end{matrix} \right)\right)=e_{-\theta}z, \end{gather} where
$\{e_{-\theta},\dot {h}_{\theta},e_{\theta}\}$ satisfy the $
sl(2,\mathbb C )$-commutation relations, and $e_{\theta}$ is a
highest root for $\dot {\mathfrak g}$.

Let $\Lambda_j$ denote the fundamental dominant integral
functionals on~$\mathfrak h$.  Any linear function $\lambda$ on~$\mathfrak h$ can be written uniquely as $\lambda =\dot{\lambda
}+\lambda (h_0)\Lambda_0$, where $\dot{\lambda}$ can be identif\/ied
with a linear function on $ \dot {\mathfrak h}$. In particular
$\delta$, the sum of the fundamental dominant integral
functionals, is given by $\delta=\dot {\delta}+\dot {g}\Lambda_0$,
where $\dot {\delta}$ is the sum of the fundamental dominant
integral functionals for the f\/inite dimensional triangular
structure (\ref{6.1}).

For $\tilde {g}\in N^{-}\cdot H\cdot N^{+}\subset\tilde {L}G$,
\begin{gather*}
\tilde {g}=l\cdot ({\rm diag}\tilde {)}\cdot u,\qquad {\rm where}\quad ({\rm diag}
\tilde {)}(\tilde {g})=\prod_0^{rk\dot {\mathfrak
g}}\sigma_j(\tilde { g})^{h_j},\end{gather*} where
$\sigma_j=\sigma_{\Lambda_j}$ is the matrix coef\/f\/icient
corresponding to $\Lambda_j$. If $\tilde {g}$ projects to $g\in
N^{-}\cdot\dot {H}\cdot N^{+}\subset LG$, then because
$\sigma_0^{h_0}=\sigma_0^{c-\dot {h}_{\theta}}$ projects to
$\sigma_ 0^{-\dot {h}_{\theta}}$, we have $g=l\cdot {\rm diag}\cdot u$,
where
\begin{gather}\label{loopdiagonal}{\rm diag}(g)=\sigma_0(\tilde {g})^{-\dot {h}_{\theta}}\prod_1^{rk\dot {
\mathfrak g}}\sigma_j(\tilde {g})^{\dot {h}_j}=\prod_1^{rk\dot
{\mathfrak g}}\left (\frac {\sigma_j(\tilde {g})}{\sigma_0(\tilde
{g})^{\check {a}_j}}\right )^{\dot {h}_j}, \end{gather} and the
$\check {a}_j$ are positive integers such that $\dot
{h}_{\theta}=\sum\check {a}_j\dot {h}_j$.

If $\tilde{g} \in \tilde{L}K$, then
$\vert\sigma_j(\tilde{g})\vert$ depends only on $g$, the
projection of $\tilde{g}$ in $LK$. We will indicate this by
writing
\begin{gather}\label{diagonalnotation}\vert\sigma_j(\tilde{g})\vert=\vert\sigma_j\vert(g).\end{gather}

In this paper we will mainly deal with generic elements in
$\tilde{L}K$ having diagonal elements with trivial $T$-component.
Thus (\ref{diagonalnotation}) has the practical consequence
(important in Sections \ref{SU(2)case} and \ref{SU(2)caseII}) that
we can generally work with ordinary loops in $K$. We record this
for later reference.

\begin{lemma}\label{notationlemma} The restriction of the projection $\tilde{L}K
\to LK$ to generic elements with diagonal terms having trivial
$T$-component is injective.
\end{lemma}

\section{Type I case}\label{typeIcase}

In this section we assume that $\dot {X}$ is a type I simply
connected and irreducible symmetric space. We let $\dot {U}$
denote the universal covering of the identity component of the
group of automorphisms of $\dot {X}$, and so on, as in the
Introduction. The irreduciblity and type I conditions imply that
$\dot {\mathfrak u}$ and $\dot {\mathfrak g}$ are simple Lie
algebras.

Exactly as in the preceding section, we introduce the af\/f\/ine
analogues $\mathfrak g=\widehat {L}\dot {\mathfrak g}$ and
$\mathfrak u=\widehat {L}\dot {\mathfrak u}$ of~$\dot {\mathfrak
g}$ and its compact real form $\dot {\mathfrak u}$, respectively,
and also the corresponding groups. We will write the corresponding
Lie algebra involution as $-(\cdot)^*$, as we typically would in a
f\/inite dimensional matrix context.

Let $\Theta$ denote the involution corresponding to the pair
$(\dot {\mathfrak u},\dot {\mathfrak k})$. We extend $\Theta$
complex linearly to $\dot {\mathfrak g}$, and we use the same
symbol to denote the involution for the Lie group $\dot {G}$. We
assume that the triangular decomposition of the preceding section
is $\Theta$-stable. We extend $\Theta$ to an involution of $L\dot
{\mathfrak g}$ pointwise, and we then extend $\Theta$ to $\widehat
{L}\dot {\mathfrak g}$ by
\[\Theta (\mu d+x+\lambda c)=\mu d+\Theta (x)+\lambda c.\]
The triangular decomposition for $\widehat {L}\dot {\mathfrak g}$
is $\Theta$-stable, and $\mathfrak t_0= \mathfrak h\cap\widehat
{L}\dot {\mathfrak k}$ is maximal abelian in $\widehat {L}\dot
{\mathfrak k}$. We let $\sigma$ denote the Lie algebra involution
$-(\cdot)^{*\Theta}$, we use the same symbol for the corresponding
group involution, and we let $\mathfrak g_0=\widehat {L}
\dot{\mathfrak g}_0$ and $G_0=\widehat {L} \dot{G}_0$ denote the
corresponding real forms.

We have def\/ined the various objects in the diagram
(\ref{groupdiagram}). The Lie algebra analogue of the diagram
(\ref{groupdiagram}) is given by
\begin{gather*}
\begin{matrix} &&\widehat {L}\dot {\mathfrak g}=\widehat {L}\dot {\mathfrak u}\oplus i\widehat {L}\dot {\mathfrak u}&&\\
&\nearrow&&\nwarrow&&\\
\widehat {L}\dot {\mathfrak g}_0=\widehat {L}\dot {\mathfrak
k}\oplus L\dot {\mathfrak p}&&&&\widehat {L}\dot {\mathfrak u}
=\widehat {L}\dot {\mathfrak k}\oplus iL\dot {\mathfrak p}\\
&\nwarrow&&\nearrow\\
&&\widehat {L}\dot {\mathfrak k}\end{matrix} \end{gather*} where
$\widehat {L}\dot {\mathfrak k}=i\mathbb R d\propto\tilde {L}\dot
{\mathfrak k}$ and $\tilde {L}\dot {\mathfrak k}=L\dot {\mathfrak
k}\oplus i\mathbb R c$. The sums in the diagram represent Cartan
decompositions. In analogy with \cite{CP}, we will write
$\mathfrak p=L \dot {\mathfrak p}$, $\mathfrak h_0=\mathfrak h
\cap \mathfrak g_0=\mathfrak t_0 \oplus \mathfrak a_0$ (relative
to the Cartan decomposition for $\mathfrak g_0$), and $\mathfrak
t=\mathfrak h \cap \mathfrak u=\mathfrak t_0 \oplus i\mathfrak
a_0$.

Our next task is to explain the diagram (\ref{Cartandiagram}).
There are isomorphisms induced by natural maps
\begin{gather}\label{compactmaps}\widehat {L}\dot {U}/\widehat {L}\dot {K}\to\tilde
{L}\dot {U}/\tilde {L}\dot {K}\to L\dot {U}/L\dot {K} \to L\dot
{X},\end{gather}and
\begin{gather}\label{noncompactmaps}\widehat {L}\dot
{G}_0/\widehat {L}\dot {K}\to\tilde {L}\dot {G}_0/\tilde {L}\dot
{K}\to L\dot {G}_0/L\dot {K}\to L\dot {X}_0.\end{gather} In each
case the f\/irst two maps are obviously isomorphisms. In the f\/irst
and second cases the third map is an isomorphism because $\dot
{X}$ and $\dot {X}_0$ are simply connected, respectively.

We will take full advantage of these isomorphisms, and
consequently there will be times when we want to use the quotient
involving hats, or tildes, and times when we want to use the
quotient not involving hats, or tildes. To distinguish when we are
using hats, we will write our group elements with hats, and
similarly with tildes. Thus $\widehat {g}$ will typically denote
an element of~$\widehat {L} \dot {G}$, whereas $g$ will typically
denote an element of~$L \dot {G}$, and unless stated otherwise,
these two elements will be related by projection.

For the natural maps \begin{gather}\label{maps1}\widehat {L}\dot
{G}/\widehat {L}\dot {G}_0\to\tilde {L}\dot {G}/\tilde {L}\dot
{G}_0\to L\dot {G}/L\dot {G}_0\to L(\dot {G}/\dot
{G}_0)\end{gather} and \begin{gather*}
\widehat
{L}\dot {G}/\widehat {L}\dot {U}\to\tilde {L}\dot {G}/\tilde
{L}\dot {U}\to L\dot {G}/L\dot {U}\to L(\dot {G}/\dot
{U})\end{gather*} in each case the third map is an isomorphism,
but the f\/irst two maps fail to be isomorphisms. For example in
(\ref{maps1}) the second map is surjective, but there is a
nontrivial f\/iber $\exp(\mathbb R c) \tilde {L}\dot {G}_0$ over the
basepoint (represented by $1$). This is the reason for the
appearance of $\tilde {L}\dot {G}/\tilde {L}\dot {G}_0$, rather
than $L(\dot {G}/\dot {G}_0)$, in the diagram
(\ref{Cartandiagram}).

There is an Iwasawa decomposition for $\widehat {L}\dot {G}$ (see
Chapter 8 of \cite{PS}), which we write as
\begin{gather}\label{Iwasawa}\widehat {L}\dot {G}\simeq N^{-}\times A \times\widehat {L}\dot {U}: \ \widehat {g}=\mathbf l
(\widehat {g})\mathbf a(\widehat {g})\mathbf u(\widehat
{g}),\end{gather} where $A=\exp(\mathfrak h_{\mathbb R})$. In
analogy with \cite{CP}, we also write $\mathbf a=\mathbf a_0
\mathbf a_1$, relative to $\exp(\mathfrak h_{\mathbb
R})=\exp(\mathfrak a_0)\exp(i\mathfrak t_0)$. There is an induced
right action \begin{gather}\label{rightaction}
\widehat {L}\dot
{U}\times (\widehat {T}\times\widehat {L}\dot {G}_0)\to\widehat
{L}\dot {U}: \ (\widehat { u},\widehat {t},\widehat {g}_0)\to\widehat
{t}^{-1}\mathbf u(\widehat {u}\widehat {g}_0 )\end{gather}
arising from the identif\/ication of $\widehat {L}\dot {U}$ with
$N^{-}A\backslash \widehat {L}\dot {G}$. We also write $A_0=A \cap
G_0$.

The Cartan embedding for the unitary type symmetric space is given
by
\[\phi : \ L(\dot {U}/\dot {K})\to\tilde {L}\dot {U}\subset\widehat {L}\dot {U}: \ \tilde {u}\tilde {L}
\dot {K}\to\tilde {u}\tilde {u}^{-\Theta},\] where we are using
the isomorphism (\ref{compactmaps}) in an essential way to express
this mapping. There is a corresponding embedding $\psi$ in the
dual case. More generally \[\phi : \ \tilde {L}\dot {G}/\tilde
{L}\dot {G}_0\to\tilde {L}\dot {G}\subset\widehat {L}\dot
{G}: \ \tilde {g}\tilde {L} \dot {G}_0\to\tilde {g}\tilde
{g}^{*\Theta},\] and the extension of $\psi$ is similarly def\/ined.

This explains the diagram (\ref{Cartandiagram}). We should note
that in what follows, in place of (\ref{groupdiagram}) and~(\ref{Cartandiagram}), and the Kac--Moody triangular decomposition
(\ref{looptriangulardecomposition}) for $\mathfrak g$, we could
simply consider the ordinary loop functor of the diagrams (0.1)
and (0.2) of \cite{CP}, and the analogue of the triangular
decomposition for $L\dot {\mathfrak g}$. But in the process we
would miss out on the interesting applications (such as Theorem
\ref{hankeltheorem}), and in analyzing the resulting Hamiltonian
systems we would inevitably be led to this Kac--Moody extended
point of view.

We are now in a position to repeat verbatim the arguments in
Sections 2--4 of \cite{CP}, supplemented with remarks concerning
Poisson structures in inf\/inite dimensions. We will summarize the
main points.

\begin{proposition} Relative to the extended real form ${\rm Im}\langle\cdot
,\cdot\rangle$ on $\mathfrak g=\widehat {L}\dot {\mathfrak g}$,
\[(\mathfrak g,\mathfrak u,\mathfrak h_{\mathbb R}\oplus
\mathfrak n^{-}) \qquad \text{and}\qquad  (\mathfrak g,\mathfrak g_0,\mathfrak
t\oplus \mathfrak n^-)\] are Manin triples, extending the finite
dimensional Manin triples $(\dot {\mathfrak g},\dot {\mathfrak
u},\dot {\mathfrak h}_{\mathbb R}\oplus \dot {\mathfrak n}^{-})$
and $(\dot {\mathfrak g},\dot {\mathfrak g}_0,\dot {\mathfrak
t}\oplus \dot {\mathfrak n}^{-})$, respectively.
\end{proposition}

We next apply the Evens--Lu construction to obtain global Poisson
structures $\Pi_X$ and $\Pi_{X_0}$ on the loop spaces $X=L\dot
{X}$ and $X_0=L\dot {X}_0$, respectively, using the isomorphisms~(\ref{compactmaps}) and~(\ref{noncompactmaps}). These Poisson
structures are given by the same formulas as in the f\/inite
dimensional cases: see~(3.1) and~(4.1) of \cite{CP}. As in the
f\/inite dimensional case, we have used the Ad-invariant symmetric
form on $\widehat {L} \dot {\mathfrak g}$ to identify $\mathfrak
p$ with a subspace of its dual (note the form is def\/inite on~$\mathfrak p$). However, in this inf\/inite dimensional context, the
inclusion $\mathfrak p \to \mathfrak p^*$ is proper, so that this
Poisson structure must be understood in a weak sense. Consequently
it is not clear that we can appeal to any general theory (e.g.\ as
in \cite{OR}) for the existence of a symplectic foliation, etc.

As in \cite{CP}, the Hilbert transform $\mathcal H\colon \mathfrak
g\to\mathfrak g$ associated to the triangular decomposition of
$\mathfrak g$ is given by
\begin{gather*}
x=x_{-}+x_0+x_{+} \mapsto \mathcal H(x)=-ix_{-}+ix_{+}.
\end{gather*}

In the following statement, we can, and do, view $\mathbf a_0$
(def\/ined following (\ref{Iwasawa})) as a function on $X=G_0/K$.

\begin{theorem} \qquad {}

{\rm (a)} The Poisson structure $\Pi_{X_0}$ has a regular
symplectic foliation (by weak symplectic manifolds), given by the
level sets of the function $\mathbf a_0$.

{\rm (b)} The horizontal parameterization for the symplectic leaf
through the basepoint is given by the map $s\colon A_0\backslash
G_0/K\to G_0/K$
\begin{gather*}
A_0 g_0 K \to s(A_0g_0K)=\mathbf a_0^{-1}g_0K,
\end{gather*}
where $g_0=\mathbf l\mathbf a_0\mathbf a_1\mathbf u$.

{\rm (c)} If we identify $T(G_0/K)$ with $G_0\times_{K}\mathfrak p$ in
the usual way, then
\begin{gather}\label{defofomega}
\omega_{1} ([g_0,x]\wedge [g_0,y])=\langle {\rm Ad}\big(\mathbf
u(g_0)^{-1}\big)\circ \mathcal H\circ {\rm Ad}(\mathbf u(g_0))(x),y\rangle
\end{gather} is a well-defined two-form on $G_0/K$.

{\rm (d)} Along the symplectic leaves, $\Pi_{X_0}^{-1}$ agrees with the
restriction of the closed two-form $\omega_{1}$.
\end{theorem}

Note that the facts that the form $\omega_1$ is closed and
nondegenerate (on the double coset space $A_0\backslash G_0/K\to
G_0/K$) is proven directly in Section 1 of \cite{CP}.

\begin{theorem} \qquad {}

{\rm (a)} The Poisson structure $\Pi_X$ has a symplectic foliation
(by weak symplectic manifolds). The symplectic leaves are
identical to the projections of the $L\dot {G}_0$-orbits, for
$L\dot {G}_0$ acting on $L\dot {U}$ as in~\eqref{rightaction}, to
$L(\dot {U}/ \dot {K})$. Let $S(1)$ denote the symplectic leaf
containing the identity.

{\rm (b)} The action of $\widehat {T}_0={\rm Rot}(S^1)\times T_0\times
\exp(i\mathbb R c)$ on $\Sigma^{\phi (L(\dot {U}/\dot {K}))}_1$ is
Hamiltonian with momentum mapping
\[\Sigma_1^{\phi (L(\dot {U}/\dot {K}))}\to (\widehat {\mathfrak t}_0)^{*}:\tilde {u}\to
\langle -\frac i2\log(a_{\phi}(\tilde {u}),\cdot\rangle ,\] where
$\tilde {u}$ has the unique triangular decomposition $\tilde
{u}=lm\tilde {a}_{\phi}l^{*\Theta}$.

{\rm (c)} The map $\Tilde{\mathbf u}\colon G_0 \to U$
\[
g_0\mapsto \mathbf u(g_0),
\]
where $\mathbf u$ is defined by \eqref{Iwasawa}, is equivariant
for the right actions of $K$ on $G_0$ and $U$, invariant under the
left action of $A_0$ on $G_0$ and descends to a $T_0$-equivariant
diffeomorphism
\begin{gather*}
\Tilde{\mathbf u}\colon  A_0\backslash G_0/ K  \to  S(1).
\end{gather*}
This induces an isomorphism of $T_0$-Hamiltonian spaces
\[
(A_0\backslash G_0/K,\omega_{1})\to \big(S(1),\Pi_X^{-1}\big),
\]
where $\omega_{1}$ is as in \eqref{defofomega}.
\end{theorem}

The symplectic foliation in part (a) can be described in a
completely explicit way in terms of triangular factorization and
the Cartan embedding $\phi$ (see \cite{C} for the f\/inite
dimensional case; the arguments there extend directly).

Throughout this paper we will focus on the generic system $S(1)$
in part (c). As we mentioned in the Introduction, the main
application which we envision is to use this Hamiltonian system to
generate useful integral formulas. In this loop context these
integrals are inf\/inite dimensional, and more infrastructure and
analysis are required to properly formulate and justify them (see~\cite{pickrell1}, especially Section~7). Even in f\/inite
dimensions, it is not known whether these type I systems have any
integrability properties (in sharp contrast to the type II case).

\section[The $S^2$ case]{The $\boldsymbol{S^2}$ case}\label{S^2case}

In this section we will do some illustrative calculations in the
simplest Type I case
\[\begin{matrix} &&G=\hat {L}SL(2,\mathbb C)\\
&\nearrow&&\nwarrow\\
G_0=\hat {L}SU(1,1)&&&&U=\hat {L}SU(2)\\
&\nwarrow&&\nearrow\\
&&K=\hat {L}U(1)\end{matrix} \] If we identify $\dot {X}_0$ with
$\Delta$ (the unit disk) and $\dot { X}$ with $\hat {\mathbb C}$
in the usual way, then from the preceding section we have maps
\begin{gather}\label{maps}L\Delta\overset {\tilde {\mathbf u}}\to L\hat {\mathbb C}\overset {
\phi}\to \tilde {L}SU(2),\end{gather} where the map $\tilde
{\mathbf u}$ is covered by the map
\[\hat {L}SU(1,1)\overset {\mathbf u}\to \hat {L}SU(2)\]
induced by the Iwasawa decomposition $g_0=\mathbf l(g_0)\mathbf
a(g_0)\mathbf u(g_0)$.

To orient the reader, we recall the nonloop case:
\[\begin{matrix} SU(1,1)&\overset {\mathbf u}\to &SU(2)\\
\downarrow&&\downarrow\\
\Delta&\to&\hat {\mathbb C}&\overset {\phi}\to &SU(2)\end{matrix}
\]
\[\begin{matrix} g_0=\frac 1{(1-Z\bar {Z})^{1/2}}\left(\begin{matrix} 1&\bar {Z}\\
Z&1\end{matrix} \right)&\to&\mathbf u(g_0)=\frac 1{(1+Z\bar
{Z})^{1/2}}\left
(\begin{matrix} 1&\bar {Z}\\
-Z&1\end{matrix} \right)\\
\downarrow&&\downarrow\\
Z&\to&-Z&\to&\frac 1{1+\vert Z\vert^2}\left(\begin{matrix} 1-\vert
Z\vert^
2&2\bar {Z}\\
-2Z&1-\vert Z\vert^2\end{matrix} \right)\end{matrix}
\] In this context $\mathbf u$ is obtained by a
Gram--Schmidt process from the rows of $g_0$, and
\[\mathbf a=\left(\frac {1+Z\bar {Z}}{1-Z\bar {Z}}\right)^{\frac 12
h_1}.\]

To calculate the symplectic form note that
\begin{gather*}
\left[g_0,X=\left(\begin{matrix} 0&\bar {x}\\
x&0\end{matrix} \right)\right]\to\frac d{dt}\Big\vert_{t=0}Z\big(g_0e^{tX}\big)\\
\qquad{}{} =\frac d{dt}\Big\vert_{t=0}\big(Z\,{\rm ch}\,(tx)+{\rm sh}\,(tx)\big)\big({\rm ch}\,(tx)+\bar {Z}\,{\rm sh}\,(tx)\big)^{
-1}=(1-Z\bar {Z})x.
\end{gather*}
 Thus a variation $\dot {Z}$ of $Z$ will
correspond to $[g_0,X]$ with $x=(1-Z\bar {Z})^{-1}\dot {Z}$. Thus
\begin{gather*}
\omega ([g_0,X]\wedge [g_0,Y])\\
\qquad{}=\omega \left(\left[\frac 1{(1-Z\bar {Z})^{1/2}}\left(\begin{matrix} 1&\bar {Z}\\
Z&1\end{matrix} \right),\left(\begin{matrix} 0&\bar {x}\\
x&0\end{matrix} \right)\right]\wedge \left[\frac 1{(1-Z\bar
{Z})^{1/2}}\left(\begin{matrix}
1&\bar {Z}\\
Z&1\end{matrix} \right),\left(\begin{matrix} 0&\bar {y}\\
y&0\end{matrix} \right)\right]\right)\\
\qquad {}=\langle \mathcal H({\rm Ad}\left(\frac 1{(1+Z\bar {Z})^{1/2}}\left(\begin{matrix} 1&
\bar {Z}\\
-Z&1\end{matrix} \right)\right)\left(\left(\begin{matrix} 0&\bar {x}\\
x&0\end{matrix} \right)\right)\\
\qquad \quad{}\wedge {\rm Ad}\left(\frac 1{(1+Z\bar
{Z})^{1/2}}\left(\begin{matrix} 1&\bar {Z}\\
-Z&1\end{matrix} \right)\right)\left(\left(\begin{matrix} 0&\bar {y}\\
y&0\end{matrix} \right)\right)\rangle=\frac i{(1-\vert Z\vert^4)}\big(\bar{\dot {Z}}Z'-\dot {Z}\bar {Z}'
\big).
\end{gather*}

Thus
\[\omega =\frac i{(1-\vert Z\vert^4)}dZ\wedge d\bar {Z}.\]

Returning to the loop case, we denote the maps in (\ref{maps}) by
\[
f(\theta )\to F(\theta )\to\phi (F)(\theta ).
\]
We have written the argument as $\theta$, as a reminder that these
are functions on $S^1$. To calculate the map $f\to F$, we need to
f\/ind the Iwasawa decomposition
\[g_0(\theta )=\frac 1{(1-f(\theta )\bar {f}(\theta ))^{1/2}}\left
(\begin{matrix} 1&\bar {f}(\theta )\\
f(\theta )&1\end{matrix} \right)=\mathbf l(z)\mathbf a\mathbf
u(\theta ),\] and remember that $\mathbf l(z)$ extends to a
holomorphic function in the exterior of $S^1$. In turn
\[\phi (F)(\theta )=\mathbf u\mathbf u^{*\Theta}=l(z)au(z),\]
where $l=\mathbf a^{-1}\mathbf l^{-1}\mathbf a$, $a=\mathbf
a^{-2}=\vert\sigma_ 0\vert^{h_0}\vert\sigma_1\vert^{h_1}$,
$u=l^{*\Theta}$. The image of $\phi (F)$ in $LSU(2)$ has the form
\[\left(\begin{matrix} \alpha (\theta )&\beta (\theta )\\
-\bar{\beta }(\theta )&\bar{\alpha }(\theta )\end{matrix}
\right)=l \big(z^{-1}\big)\left(\frac
{\vert\sigma_1\vert}{\vert\sigma_0\vert}\right)^{h_1}u(z )\] (see
(\ref{loopdiagonal})).

The Iwasawa decomposition of $g_0$ (the special self-adjoint
representative above) is equivalent to
\[g_0^{*}\mathbf l(g_0)^{-*}\mathbf a(g_0)^{-1}=g_0\mathbf l(g_0)^{-*}\mathbf a
(g_0)^{-1}=\mathbf u(g_0).
\]
Write $\mathbf l^{-*}=\left(\begin{matrix} a&b\\
c&d\end{matrix} \right)$, so that $a$, $b$, $c$, $d$ are holomorphic
functions in $D$, $a(0)=d(0)=1$, and $c(0)=0$. Then
\[\left(\begin{matrix} 1&\bar {f}\\
f&1\end{matrix} \right)\left(\begin{matrix} a&b\\
c&d\end{matrix} \right)\left(\begin{matrix} a_0^{-1}&0\\
0&a_0\end{matrix} \right)\]
is of the form $\left(\begin{matrix} A&B\\
-\bar {B}&\bar {A}\end{matrix} \right)$. This implies
\begin{gather}\label{equations}a+\bar {f}c=\big(\bar {f}\bar {b}+\bar {d}\big)a_0^2, \qquad
fa+c=-\big(\bar {b}+f\bar {d}\big)a_0^2.\end{gather} As a reminder,
these are equations for functions def\/ined on $S^1$.

Let $\mathcal H_0=P_{+}-P_{-}$, where for a scalar function
$g=\sum g_ nz^n$, $P_{+}g=\sum\limits_{n\ge 0}g_nz^n$. We take the
conjugate of the f\/irst equation in (\ref{equations}) and rewrite
it as
\[
-\mathcal H_0\big(\bar {a}+da_0^2\big)+2=f\mathcal H_0\big(ba_0^2+\bar {c}\big).
\]
This is equivalent to
\[
-\big(\bar {a}+da_0^2\big)+2=\mathcal H_0f\mathcal H_0\big(ba_0^2+\bar {c}\big).
\]
The second equation in (\ref{equations}) is equivalent to
\[
\bar {f}\big(\bar {a}+da_0^2\big)=-\big(ba_0^2+\bar {c}\big).
\]
These two equations imply
\[
\big(ba_0^2+\bar {c}\big)=-2(1-\bar {f}\mathcal H_0f\mathcal H_0)^{-1}(\bar {f}
)=-2\bar {f}(1-\mathcal H_0f\mathcal H_0\bar {f})^{-1}(1).
\] Note
that the inverse on the right exists, because $\sup\{\vert
f(z)\vert :z\in S^1\}<1$. This determines $b a_0^2$ and $c$, by
applying $P_{\pm}$.

We now see that
\[\bar {a}+da_0^2=2\big(1+\mathcal H_0f\mathcal H_0\bar {f}(1-\mathcal H_0f\mathcal H_
0\bar {f})^{-1}(1)\big)=2(1-\mathcal H_0f\mathcal H_0\bar
{f})^{-1}(1).
\] Note that $1+a_0^2$ is the zero mode of the right
hand side, so that in principle we have determi\-ned~$a_0$, and
$\mathbf l$. This form of the solution does not explain in a clear
way why the zero mode of the right hand side is $>2$.

To summarize, let
\[
h=2(1-\mathcal H_0f\mathcal H_0\bar {f})^{-1}(1)
\]
(this is a well-def\/ined function on $S^1$, and we do not know much
more about it). Then
\[
\bar {c}=-P_{-}(\bar {f}h), \qquad
ba_0^2=-P_{+}(\bar {f}h), \qquad \bar {a}-1=P_{-}(h), \qquad
1+da_0^2=P_{+}(h).
\] This implies the following

\begin{proposition}
\[\mathbf u=g_0\mathbf l^{-*}\mathbf a^{-1} =(1-f\bar {f})^{-1/2}a_0^{-1}\left(\begin{matrix} (1+P_{-}h-fP_{-}(\bar {
f}h))^{*}&-(\bar {f}+\bar {f}P_{-}h-P_{-}(\bar {f}h))\vspace{1mm}\\
(\bar {f}+\bar {f}P_{-}h-P_{-}(\bar
{f}h))^{*}&1+P_{-}h-fP_{-}(\bar { f}h)\end{matrix} \right)\] and
\[
F=\left(\frac {\bar {f}+\bar {f}P_{-}h-P_{-}(\bar {f}h)}{1+P_{-}h-fP_{
-}(\bar {f}h)}\right)^{*}\in L\hat {\mathbb C}.
\]
\end{proposition}

These general formulas are not especially enlightening. However,
the $SU(2)$ case considered below suggests that there might be
some special cases of these formulas which are tractable.

\section{Type II case}\label{typeIIcase}

In the type II case there is more than one reasonable
interpretation of the diagram (\ref{groupdiagram}). The
dif\/ferences between the possibilities are minor, but potentially
confusing. We will brief\/ly describe a f\/irst possibility, which
leads to diagram (\ref{Cartandiagram}), but we will then consider
a second possibility, which is more elementary in a technical
sense, and we will pursue this in detail.

Throughout this section $\dot {K}$ denotes a simply connected
compact Lie group with simple Lie algebra $\dot {\mathfrak k}$,
$\dot {X}=\dot {K}$, viewed as a symmetric space, $\dot {U}=\dot
{K} \times \dot {K}$, and $\dot {\mathfrak g}=\dot {\mathfrak
k}^{\mathbb C}\oplus\dot {\mathfrak k}^{\mathbb C}$.

In the f\/irst interpretation of diagram (\ref{groupdiagram}),
$\mathfrak g=\hat {L}\dot {\mathfrak g}$ is def\/ined in the
following way. We f\/irst def\/ine a central extension
\[0\to \mathbb C c \to \tilde {L}\dot {\mathfrak g}\to L\dot {\mathfrak g}\to 0.\]
As a vector space
\[\tilde {L}\dot {\mathfrak g}=L\dot {\mathfrak g}\oplus \mathbb C c;\]
the bracket is def\/ined as in (\ref{bracket}), where the form
$\langle\cdot , \cdot\rangle$ is the sum of the normalized
invariant symmetric forms for the two $\dot {\mathfrak k}^{\mathbb
C}$ factors:
\[\langle (x,y),(X,Y)\rangle =\langle x,y\rangle +\langle X,Y\rangle
.\] Then $\hat {L}\dot {\mathfrak g}=\mathbb C d\propto\tilde
{L}\dot {\mathfrak g}$, and
\[\Theta (\lambda d+(x,y)+\mu c)=\lambda d+(y,x)+\mu c.\]

The Lie algebra analogue of diagram (\ref{groupdiagram}) is
\[\begin{matrix} &&\mathfrak g=\hat {L}\dot {\mathfrak g}\\
&\nearrow&&\nwarrow\\
\hat {L}\dot {\mathfrak k}\oplus \{(x,-x):x\in iL\dot {\mathfrak
k}\}&&&&\hat {L}\dot {\mathfrak k}
\oplus \{(x,-x):x\in L\dot {\mathfrak k}\}\\
&\nwarrow&&\nearrow\\
&&\hat {L}\dot {\mathfrak k}\end{matrix} \] where $\hat {L}\dot
{\mathfrak k}=i\mathbb R d\propto\tilde {L}\dot {\mathfrak k}$ and
$\tilde { L}\dot {\mathfrak k}=\{(x,x):x\in L\dot {\mathfrak
k}\}\oplus i\mathbb R c$.

At the group level $G=\hat {L}\dot {G}=\mathbb C \propto\tilde
{L}\dot {G}$ where $ \tilde {L}\dot {G}$ is an extension of $L\dot
{G}$ by $\mathbb C^{*}$; precisely, $\tilde {L}\dot {G}$ is a
quotient
\[
0\to \mathbb C^{*}\to\tilde {L}\dot {K}^{\mathbb C}\times\tilde {L}\dot {K}^{\mathbb C}
\to\tilde {L}\dot {G}\to 0,
\] where $\lambda\in \mathbb C^{*}$ maps
antidiagonally, $\lambda\to (\lambda^c,\lambda^{-c})$.

As in the type I case, there are isomorphisms
\[\hat {L}\dot {U}/\hat {L}\dot {K}\to\tilde {L}\dot {U}/\tilde {L}\dot {K}\to L\dot {U}/L\dot {K}\to L\dot {K},\]
where the last map is given by $(k_1,k_2)\to k_1k_2^{-1}$.
The Cartan embedding is given by
\[\phi : \ X=L\dot {K}\to\tilde {L}\dot {U}\subset\hat {L}\dot {U}:\ k\to\widetilde{(k_1,k_2)}\widetilde{
(k_1,k_2)}^{-\Theta},
\] where $k=k_1k_2^{-1}$. The dual map $\psi$
is described in a similar way. This leads to the diagram
(\ref{Cartandiagram}) in this type II case.

This f\/irst interpretation of diagram (\ref{groupdiagram}) is
somewhat inconvenient, because unlike the f\/inite dimensional case,
$X \ne K$, and $U \ne K \times K$. In the remainder of this paper
we will consider a setup where these equalities do hold. It will
be easier to compare this setup with the f\/inite dimensional case.
The modest price we pay is that, in this second interpretation,
$X$ is a covering of $L\dot {K}$ (also, as a symmetric space, the
invariant geometric structure is of Minkowski type, rather than
Riemannian type, but this geometric structure is irrelevant for
our purposes).

From now on, we set $K=\widehat {L} \dot{K}$, as in Section
\ref{loopgroups}. We henceforth understand the diagram
(\ref{groupdiagram}) to be
\[\begin{matrix} &&&G=K^{\mathbb C}\times K^{\mathbb C}\\
&\nearrow&&&\nwarrow\\
G_0=\{g_0=(g,g^{-*}):g\in K^{\mathbb C}\}&&&&&U=K\times K\\
&\nwarrow&&&\nearrow\\
&&&\Delta (K)=\{(k,k): k \in K\}\end{matrix}\] where $\Delta
(K)=\{(k,k)\colon k\in K\}$, $G_0=\{g_0=(g,g^{-*})\colon g\in
K^\mathbb C\}$, $K=\widehat {L} \dot {K}$, $G=\widehat {L} \dot
{G}$, and the involution $\Theta$ is the outer automorphism
$\Theta ((g_1,g_2))=(g_2,g_1)$. Also \[ X_0=G_0/\Delta (K)\simeq
K^{\mathbb C}/K,\text{ and } X=U/\Delta (K)\simeq K,
\]
where the latter isometry is $(k_1,k_2)\Delta (K)\mapsto
k=k_1k_2^{-1}$. As in \cite{CP} we will use superchecks to
distinguish structures for $k^{\mathbb C}$ versus those for
$\mathfrak g$.

We f\/ix a triangular decomposition
\begin{gather}\label{groupcase1}
\check {\mathfrak g}=\mathfrak k^{\mathbb C}=\check{\mathfrak
n}_{-} + \check{\mathfrak h} + \check {\mathfrak n}_{+}.
\end{gather}
This induces a $\Theta$-stable triangular decomposition for
$\mathfrak g$
\begin{gather}\label{groupcase2}
\mathfrak g = \underbrace{(\check{\mathfrak
n}^-\times\check{\mathfrak n}^-)}_{\mathfrak
n^-}+\underbrace{(\check{\mathfrak h}\times\check{\mathfrak
h})}_{\mathfrak h}+\underbrace{(\check{\mathfrak
n}^+\times\check{\mathfrak n}^+)}_{\mathfrak n^+}.
\end{gather}
Let $\check {\mathfrak a}=\check {\mathfrak h}_{\mathbb R}$ and
$\check {\mathfrak t}=i\check {\mathfrak a}$. Then
\[ \mathfrak t_0=\{(x,x):x\in\check {\mathfrak
t}\},\qquad\text{\rm and}\qquad \mathfrak a_0=\{(y,-y):y\in\check
{\mathfrak a}\}.
\]
The standard Poisson Lie group structure on $U=K\times K$ induced
by the decomposition in (\ref{groupcase2}) is then the product
Poisson Lie group structure for the standard Poisson Lie group
structure on~$K$ induced by the decomposition (\ref{groupcase1}).

Let us denote the Poisson Lie group structure on $K$ by $\pi_K$
and the Evens--Lu homogeneous Poisson structure on $X=K$ by
$\Pi_X$. The formal identif\/ication of $\mathfrak k$ with its dual
via the invariant form allows us to view the Hilbert transform
$\check{\mathcal H}$ associated to (\ref{groupcase1}) as an
element of $\mathfrak k\wedge\mathfrak k$. As a~bivector f\/ield
\[ \pi_K = \check{\mathcal H}^r-\check{\mathcal H}^l,
\]
where $\check{\mathcal H}^r$ (resp. $\check{\mathcal H}^l$)
denotes the right (resp. left) invariant bivector f\/ield on $K$
generated by~$\check{H}$, whereas $\Pi_K=\check{\mathcal
H}^r+\check{\mathcal H}^l$.

Just as in the Type I case, the arguments of Sections 2--4 of
\cite{CP} apply verbatim. We will focus on the new issues which
arise.

As we pointed out in the Introduction, the f\/irst thing to note is
that Theorem 5.1 of \cite{CP} does not hold in this context. The
symplectic leaves for the Poisson Lie group structure on $K$ are
f\/inite dimensional, whereas the symplectic leaves for the Evens--Lu
Poisson structure are f\/inite codimensional. Thus these structures
are fundamentally dif\/ferent.

As in \cite{CP}, we will write
\begin{gather}\label{su2notation} k(\zeta
)=\left(\begin{matrix}
1&0\\
\zeta&1\end{matrix} \right)\left(\begin{matrix} a(\zeta)&0\\
0&a(\zeta)^{-1}\end{matrix} \right)\left(\begin{matrix} 1&-\bar{\zeta}\\
0&1\end{matrix} \right),\end{gather} where
$a(\zeta)=(1+\vert\zeta\vert^2)^{-1/2}$. Given a simple positive
root $\gamma$, $i_{\gamma}:SU(2) \mapsto K$ denotes the root
subgroup inclusion (as in (\ref{roothom})), and
\[r_{\gamma}=i_{\gamma}\left(\left(\begin{matrix} 0&i\\i&0\end{matrix}\right)\right),\]
a f\/ixed representative for the corresponding Weyl group
ref\/lection.

\begin{conjecture}\label{Lutheorem} Fix $w\in W$.

{\rm (a)}  The submanifold $\check {N}^{-}\cap w^{-1}\check
{N}^{+}w\subset \check {N}^{-}$ is $\check {T}$-invariant and
symplectic.

Fix a representative $\mathbf w$ for $w$ with minimal
factorization $\mathbf w=r_n\cdots r_1$, in terms of simple reflections
$r_j=r_{\gamma_j}$ corresponding to simple positive roots
$\gamma_j$.  Let $\mathbf w_j=r_j\cdots r_1$.

{\rm (b)}  The map
\[\mathbb C^n\to N^{-}\cap w^{-1}N^{
+}w: \ \zeta =(\zeta_n,\dots,\zeta_1)\to l(\zeta),
\] where
\[\mathbf w_{n-1}^{-1}i_{\gamma_
n}(k(\zeta_n))\mathbf w_{n-1}\cdots \mathbf
w_1^{-1}i_{\gamma_2}(k(\zeta_2))\mathbf
w_1i_{\gamma_1}(k(\zeta_1))=l(\zeta )au\] is a diffeomorphism.

{\rm (c)} In these coordinates the restriction of $\omega$ is given by
\begin{gather}\label{Luform}\omega\vert_{N^{-}\cap w^{-1}N^{+}w}=\sum_{j=1}^n\frac
i{\langle\gamma_ j,\gamma_j\rangle}\frac
1{(1+\vert\zeta_j\vert^2)}d\zeta_j\wedge
d\bar{\zeta}_j,\end{gather} the momentum map is the restriction
of $-\langle\frac i2\log(a ),\cdot\rangle$, where \[a(k(\zeta
))=\prod_{j=1}^n \big(1+\vert\zeta_j\vert^2\big)^{-\frac 12 w_{j-1}^{-1}
h_{\gamma_j}w_{j-1}},\]
 and Haar measure (unique up to a constant) is given by
\begin{gather*}
d\lambda_{N^{-}\cap w^{-1}N^{+}w}(l)=\prod_{j=1}^n
\big(1+\vert\zeta_j\vert^2\big)^{\check{\delta}(w_{j-1}^{-
1}h_{\gamma_j}w_{j-1})-1}d\lambda(\zeta_j)\\
\phantom{d\lambda_{N^{-}\cap w^{-1}N^{+}w}(l)}{} =\prod_{1\le i<j\le
n}\big(1+\vert\zeta_j\vert^2\big)^{-\gamma_i(w_{i-1}
w_{j-1}^{-1}h_{\gamma_j}w_{j-1}w_{i-1}^{-1})}d\lambda(\zeta_j).
\end{gather*}
where $\check{\delta }=\sum\check{\Lambda}_j$, the sum of the
dominant integral functionals for $\check {g}$, relative to~\eqref{groupcase1}.

{\rm (d)} Let $C_{\mathbf w}$ denote the symplectic leaf through
$\mathbf w$, with respect to $\pi_K$, with the negative of the
induced symplectic structure. Then left translation by $\mathbf
w^{-1}$ induces a symplectomorphism from~$C_{\mathbf w}$, with its
image in $(S(1),\omega)$, which is identified with $\check
{N}^{-}\cap w^{-1}\check {N}^{+}w \subset \check {N}^{-}$.
\end{conjecture}

\begin{proposition}\label{Luprop}The following are true:

{\rm (1)} part {\rm (b)} of the conjecture;

{\rm (2)} the formulas for $a(k(\zeta ))$ and Haar measure in part {\rm (c)};

{\rm (3)} the right hand side of \eqref{Luform} equals the image of the
symplectic structure for $C_{\mathbf w}$ with respect to the map
in part {\rm (d)}; and

{\rm (4)} the momentum maps for $\omega$ and the form in~{\rm (3)} do agree, and are given by the formula in~{\rm (c)}.
\end{proposition}

Thus the basic open question is whether (d) holds. This is known
to be true in the f\/inite dimensional case (see Theorems 5.1 and
5.2 of \cite{CP}).

\begin{proof} The proof uses a number of facts which are recalled in Appendix A.
We will freely use the notation which is used there. One technical
point which emerges is that we are currently using~$K$ to denote
an extension of the real analytic completion of $L\dot {K}$ (and
similarly for its Lie algebra, etc), whereas in the Appendix $K$
is the restriction of the extension to the polynomial loop group~$L_{\rm pol}\dot {K}$. Since all the root homomorphisms map into the
extension over the polynomial loop group, we will simply replace
$K$ by this small subgroup, rather than introducing more notation.

Via the projection,
\[K\to K/\check {T}=\check {G}/\check {B}_{+},\]
$S(1)$ is identif\/ied with $\Sigma_1$, the unique open Birkhof\/f
stratum in the f\/lag manifold.  There is a~surjective map
\begin{gather*}
SL(2,\mathbb C )\times \cdots \times SL(2,\mathbb C )\to w^{-1}\bar {C}_w:\nonumber\\
\qquad (g_n,\dots ,g   _1)\to \mathbf w_{n-1}^{-1}i_{\gamma_n}(g_n)\mathbf
w_{n-1}\cdots \mathbf w_1^{-1}i_{ \gamma_2}(g_2)\mathbf
w_1i_{\gamma_1}(g_1)\check {B}_{+}.
\end{gather*}
This map is obtained by composing the Bott--Samelson
desingularization of the Schubert va\-rie\-ty~$\bar{C}_w$ by left
translation by $\mathbf w^{-1}$, as in Corollary~\ref{dougA18}.
This map has the remarkable property that the notion of generic is
compatible with factorization:  the preimage of $\mathbf
w^{-1}C_w$ is precisely
\[SL(2,\mathbb C )'\times \cdots \times SL(2,\mathbb C )',\]
where
\[SL(2,\mathbb C )'=\left\{\left(\begin{matrix} a&b\\
c&d\end{matrix} \right):a\ne 0\right\}.\] This follows from Corollary~\ref{dougA18}.

In terms of the af\/f\/ine coordinate $\check {N}^{-}$ for $\check {G}
/\check {B}^{+}$,
\begin{gather*}
\check {N}^{-}\cap w^{-1}\check {N}^{+}w\leftrightarrow w^{-1}C_
w.
\end{gather*}
We thus obtain a surjective map
\[SL(2,\mathbb C )'\times \cdots \times SL(2,\mathbb C )'\to
\check {N}^{-}\cap w^{-1}\check {N}^{+}w: \ (g_n,\dots,g   _1)\to l,
\]
where
\[
\mathbf w_{n-1}^{-1}i_{\gamma_n}(g_n)\mathbf w_{n-1}\cdots \mathbf w_1^{-1}i_{
\gamma_2}(g_2)\mathbf w_1i_{\gamma_1}(g_1)\check
{B}^+=l\check{B}^+.
\]

By part (a) of Corollary \ref{dougA18b}, this map induces a
parameterization
\begin{gather*}
\mathbb C^n\leftrightarrow S(1)^{SU(2)} \times \cdots \times
S(1)^{SU(2)} \to \check {N}^{-}\cap w^{-1}\check {N}^{+}w,\\
(\zeta_n,\dots,\zeta_1)\leftrightarrow (k_n(\zeta_n),\dots,k_1(\zeta_
1))\to l(\zeta ),
\end{gather*} where
\[
\mathbf w_{n-1}^{-1}i_{\gamma_n}(k_n)\mathbf w_{n-1}\cdots i_{\gamma_1}(k_1)=l(\zeta
)\prod_1^na(\zeta_j)^{{\rm Ad}(w_{j-1}^{-1})(h_{\gamma_j})}u,
\]
$u\in\check {N}^{+}$. This is part (b) of the Theorem.

The formula for $a(k(\zeta))$ in part (c) follows from Proposition
\ref{dougA26}. The formula for the Haar measure follows from part~(b) of Corollary \ref{dougA18b}.

Finally part (3) of the Proposition is equivalent to Lu's
factorization result, Theorem~3.4 of~\cite{Lu}, and part (4)
follows from the results in Sections 2--4 of \cite{CP} (which, as
we have already noted, are valid in the loop context), and Theorem~3.4 of~\cite{Lu}.
\end{proof}

\section{The SU(2) case. I}\label{SU(2)case}

To understand the signif\/icance of Proposition~\ref{Luprop}, we
will now spell out its meaning in the simplest case, where $\dot
{K}=SU(2)$.

In doing explicit calculations, it is convenient to work with
ordinary loops, rather than lifts in~$\tilde{L}SU(2)$. Thus in
this section, and the next, we will identify $S(1)$ with its
projection to~$LSU(2)$ (see Lemma \ref{notationlemma}). We will
continue to denote this projection by $S(1)$.

In this case there is an (outer) automorphism of $\tilde {L}\dot
{G}$ which interchanges the simple roots $\alpha_0$ and
$\alpha_1$.  At the level of loops, this automorphism is realized
by conjugation by a multivalued loop,
\begin{gather}\label{automorphism}{\rm conj}\left(\left(\begin{matrix} 0&iz^{1/2}\\
iz^{-1/2}&0\end{matrix} \right)\right): \ \left(\begin{matrix} a&b\\
c&d\end{matrix} \right)\to\left(\begin{matrix} d&cz^{-1}\\
bz&a\end{matrix} \right).\end{gather} The root subgroup
corresponding to $\alpha_1$ is $SL(2,\mathbb C )$, the constants,
and the root subgroup corresponding to $\alpha_0$ is the image of
$SL(2,\mathbb C )$ under this automorphism (see chapter 5 of
\cite{PS}).

In \cite{PS} there is a relatively explicit realization of the
groups $\tilde {L}SU(2,\mathbb C)$ and $\tilde {L}SL(2,\mathbb C
)$.  In this approach a loop $g\in LSL(2,\mathbb C )$ is
identif\/ied with a multiplication operator on $\mathcal
H=L^2(S^1,\mathbb C^2)$. Relative to the Hardy polarization
$\mathcal H=\mathcal H_{+}\oplus \mathcal H_{-}$,
\begin{gather}\label{segalmatrixnotation}g=\left(\begin{matrix} A&B\\
C&D\end{matrix} \right), \end{gather} where $A$ (or $D$) is the
classical Toeplitz operator and $C$ (or $B$) is the classical
Hankel operator associated to $g$. The extension $\tilde
{L}SL(2,\mathbb C )\to LSL(2,\mathbb C )$ is the $\mathbb C^{
*}$-bundle associated to the pullback of the determinant line
bundle, relative to the mapping
\[LSL(2,\mathbb C )\to {\rm Fred}(\mathcal H_{+}): \ g\to A(g). \]
The holomorphic function $\sigma_0$ on $\tilde {L}SL(2,\mathbb C
)$ is, viewed as a section of a line bundle, `$\det A(\tilde {g})$'.

Suppose that $g\in LSU(2)$, and $\tilde {g}\in\tilde {L}SU(2)$ is
a lift, which is uniquely determined up to multiplication by an
element of the unitary center $\exp(i\mathbb R c)$.  Then
\[\vert\sigma_0(\tilde {g})\vert^2=\vert\sigma_0\vert^2(g)=\det A(g
)^{*}A(g)=\det(1+Z(g)^{*}Z(g))^{-1},\] where $Z=CA^{-1}$, and
\begin{gather}\label{sigma1}\vert\sigma_1\vert^2(g)=\vert\sigma_0\vert^2\left(\left(\begin{matrix} z^{1
/2}&0\\
0&z^{-1/2}\end{matrix} \right)g\left(\begin{matrix} z^{-1/2}&0\\
0&z^{1/2}\end{matrix} \right)\right). \end{gather}

The simple ref\/lections corresponding to the simple roots
$\alpha_0$ and $\alpha_1$ are represented by the group elements
\[s_0=\left(\begin{matrix} 0&iz^{-1}\\
iz&0\end{matrix} \right),\qquad {\rm and}\qquad s_1=\left(\begin{matrix} 0&i\\
i&0\end{matrix} \right),\] respectively.  We denote their images
in the Weyl group by $\bar {s}_i$.  The Weyl group (the af\/f\/ine
Weyl group of $(\dot {\mathfrak g} ,\dot {\mathfrak h})$) has the
structure
\[
W=\mathbb Z_2\bar {s}_0\propto\mathbb Z(\bar {s}_0\bar {s}_1)=\mathbb Z_
2\bar {s}_1\propto\mathbb Z(\bar {s}_0\bar {s}_1)=\dot {W}\propto
{\rm Hom}(\mathbb T,\dot {T}).
\] Minimal factorizations in the Weyl
group must simply alternate the~$\bar {s}_i$.  This leads to two
possible inf\/inite minimal sequences of simple roots, the two
possibilities depending upon whether one begins with $\alpha_0$ or~$\alpha_1$.  These are equivalent via the automorphism above. In
the following theorem we will spell out Proposition~\ref{Luprop}
for the f\/irst possibility.

\begin{theorem}\label{LSU(2)case} Let $w_1=s_0$, $w_2=s_1s_0$,
$w_3=s_0s_1s_0$, \dots.  Then for $n>0$,

{\rm (a)} \[N^{-}\cap w_n^{-1}N^{+}w_n=\left\{l=\left(\begin{matrix}
1&\sum\limits_{j=1}^nx_jz^{
-j}\\
0&1\end{matrix} \right): \ x_j\in \mathbb C\right\}.\]

{\rm (b)} For the diffeomorphisms \[\mathbb C^n\to \mathbb
C^n: \ (\zeta_1,\dots ,\zeta_n)\to x^{(n)}=\sum_{j=1}^
nx_j^{(n)}(\zeta_1,\dots ,\zeta_n)z^{-j}\] arising from the
isomorphism in part {\rm (1)} of Proposition~{\rm \ref{Luprop}}, and the
parameterization in part {\rm (a)},
\[x_j^{(n)}(\zeta_1,\dots,\zeta_n)=x_j^{(N)}(\zeta_1,\dots,\zeta_n,0,\dots
,0),\] $n<N$ (hence we will often suppress the superscript), and
\[x_j^{(n)}(\zeta_1,\dots,\zeta_n)=x_1(\zeta_j,\dots,\zeta_n,0,\dots).\]

{\rm (c)} In terms of the correspondence of $\zeta$ with $g\in LSU(2)$
and $l\in N^{-}$, arising from the isomorphism in part {\rm (1)} of
Proposition {\rm \ref{Luprop}},
\[\vert\sigma_0\vert^2(g)=\frac 1{\det\left(1+B\left(\sum\limits_1^n x_jz^j\right)B\left(\sum\limits_1^
nx_jz^j\right)^{*}\right)}=\prod_1^n\frac 1{(1+\vert\zeta_j\vert^2)^j}\] and
\[\vert\sigma_1\vert^2(g)=\frac 1{\det\left(1+B\left(\sum\limits_1^{n-1}x_{j+1}z^j\right)
B\left(\sum\limits_1^{n-1}x_{j+1}z^j\right)^{*}\right)}=\prod_1^{n-1}\frac
1{(1+\vert\zeta_{ j+1}\vert^2)^{j-1}},
\] where $B(\cdot)$ is
defined as in Theorem~{\rm \ref{hankeltheorem}}.  In particular
\[
\frac {\vert\sigma_0\vert}{\vert\sigma_1\vert}=\prod_1^n\frac 1{
(1+\vert\zeta_j\vert^2)^{1/2}}.
\]

{\rm (d)} More generally, for $0\le l<n$,
\[\det\left(1+B\left(\sum_1^{n-l}x_{j+l}z^j\right)B\left(\sum_1^{n-l}x_{j+l}z^j\right)^{*}\right)=\prod_
1^{n-l}\big(1+\vert\zeta_{j+l}\vert^2\big)^j.\]

{\rm (e)}
\[\prod_{j=1}^nd\lambda (x_j)=
\big(1+\vert\zeta_2\vert^2\big)^2\big(1+\vert\zeta_2\vert^2\big)^4\cdots \big(1+\vert\zeta_
n\vert^2\big)^{2(n-1)}\prod_{j=1}^nd\lambda (\zeta_j).\]
\end{theorem}

\begin{proof} Part (a) is a direct
calculation. If $n=2m$
\[w_n=(s_1s_0)^m=(-1)^m\left(\begin{matrix} z&0\\
0&z^{-1}\end{matrix} \right)^m.\]
Thus if $u=\left(\begin{matrix} a&b\\
c&d\end{matrix} \right)\in N^{+}$, then
\[w_n^{-1}uw_n=\left(\begin{matrix} a&z^{-2m}b\\
z^{2m}c&d\end{matrix} \right)\in N^{-}\] implies $a=d=1$, $c=0$,
and $b=\sum\limits_0^{2m-1}b_jz^j$.  This implies (a) when $n$ is even.
The odd case is similar.

Before taking on the other parts of the Theorem, we need to
understand what part (a) says in terms of the isomorphism of part
(b) of Conjecture~\ref{Lutheorem}. It is straightforward to
calculate that
\[w_{j-1}^{-1}i_{\gamma_j}(k(\zeta_j))w_{j-1}=a(\zeta_j)\left(\begin{matrix} 1&
\zeta_jz^{-j}\\
-\bar{\zeta}_jz^j&1\end{matrix} \right). \] This implies that
\[g=w_{n-1}^{-1}i_{\gamma_n}(k(\zeta_n))w_{n-1} \cdots i_{\gamma_1}(k(\zeta_
1))=a(\zeta_n)\left(\begin{matrix} 1&\zeta_nz^{-n}\\
-\bar{\zeta}_nz^n&1\end{matrix}
\right)\cdots a(\zeta_1)\left(\begin{matrix} 1&\zeta_1z^{
-1}\\
-\bar{\zeta}_1z&1\end{matrix} \right).\] If we write
\[g=a(\zeta_n)\cdots a(\zeta_1)\left(\begin{matrix} \alpha_n&\beta_n\\
\gamma_n&\delta_n\end{matrix} \right),\] then there is a recursion
relation
\begin{gather}
\left(\begin{matrix} \beta_{n+1}\\
\delta_{n+1}\end{matrix} \right)=\left(\begin{matrix} 1&\zeta_nz^{-n-1}\\
-\bar{\zeta}_nz^{n+1}&1\end{matrix} \right)\left(\begin{matrix} \beta_n\\
\delta_n\end{matrix} \right).\label{2.21}
\end{gather}

In terms of the isomorphism in part (1) of Proposition
\ref{Luprop}, part (a) implies that
\begin{gather}\label{2.22}
\left(\begin{matrix} 1&-\sum\limits_1^nx^{(n)}_jz^{-j}\\
0&1\end{matrix} \right)\left(\begin{matrix} \alpha_n&\beta_n\\
\gamma_n&\delta_n\end{matrix} \right)
=\left(\begin{matrix} \alpha_n-\left(\sum\limits_1^nx^{(n)}_jz^{-j}\right)\gamma_n&\beta_
n-\left(\sum\limits_1^nx^{(n)}_jz^{-j}\right)\delta_n\\
\gamma_n&\delta_n\end{matrix} \right)\end{gather} is an (entire)
holomorphic function of $z$.  In particular $\gamma_n$ and
$\delta_n$ must be holomorphic functions of $z$, and
\begin{gather}\label{2.23}\sum_1^nx^{(n)}_j(\zeta_1,\dots,\zeta_n)z^{-j}=\big(\delta^{-1}_n\beta_
n\big)_{-},\end{gather} where $(\cdot )_{-}$ denotes the singular
part (at $z=0$).  The holomorphicity of (\ref{2.22}) can be
checked directly as follows.  The recursion relation (\ref{2.21})
shows that $\delta_n$ is of the form $1+\sum_1^nd_jz^j$, and
$\beta_n$ is of the form $\sum_ 1^nb_jz^{-j}$. Since
$\gamma_n^{}=-\beta_n^{*}$ on $S^1$, this shows $\gamma_n$ and $
\delta_n$ are holomorphic functions of $z$.  It also shows the
$x^{(n)}_j$ are well-def\/ined by (\ref{2.23}). The relation
(\ref{2.23}) implies the $(1,2)$ entry of (\ref{2.22}) is
holomorphic. Also (\ref{2.23}) implies the $(1,1)$ entry of
(\ref{2.22}) is of the form
\begin{gather*}
\alpha_n-\delta_n^{-1}\beta_n\gamma_n+{\rm holomorphic} =\delta_n^{-1}(\alpha_n\delta_n-\beta_n\gamma_n)+{\rm hol.}\\
\qquad{}=\delta_n^{-1}(\alpha_n\alpha_n^{*}+\beta_n\beta_n^{*})+{\rm hol.} =({\rm const})\delta_n^{-1}+{\rm hol.}={\rm holomorphic}.
\end{gather*}

We now consider part (b).  We will need several Lemmas.

\begin{lemma}\label{2.25} \qquad {}

{\rm (a)} The $x^{(n)}$ satisfy the
recursion relation
\begin{gather}\label{2.26a}x^{(n+1)}=\left(\big(x^{(n)}+\zeta_{n+1}z^{-n-1}\big)\sum_{p=0}^n(\bar\zeta_{
n+1}x^{(n)}z^{n+1})^p\right)_{-}\\
\label{2.26b}
\phantom{x^{(n+1)}}{} =\left(\sum_{p=0}^n\bar\zeta_{n+1}^p\big(x^{(n)}\big)^{p+1}\big(1+\vert\zeta_{
n+1}\vert^2\big)z^{p(n+1)}\right)_{-}+\zeta_{n+1}z^{-n-1}.
\end{gather}

{\rm (b)} $x^{(n)}$ can be replaced by $x^{(n)}+h(z)$, where $h(z)$ is a
holomorphic function, without changing the recursion.
\end{lemma}

\begin{proof} The recursion relation (\ref{2.21}) and
(\ref{2.23}) imply that $x^{(n+1)}$ is the singular part of
\begin{gather}
\big(\beta_n+\delta_n\zeta_{n+1}z^{-n-1}\big)\big(\delta_n-\bar{\zeta}_{n+1}
z^{n+1}\beta_n\big)^{-1}
=\big(\beta_n\delta_n^{-1}+\zeta_{n+1}z^{-n-1}\big)\big(1-\bar{\zeta}_{n+1}
\beta_n\delta_n^{-1}z^{n+1}\big)^{-1}\nonumber\\
\label{2.27}
\qquad {}=\big(\beta_n\delta_n^{-1}+\zeta_{n+1}z^{-n-1}\big)\sum_{p=0}^{\infty}\big(
\bar{\zeta}_{n+1}x^{(n)}z^{n+1}\big)^p.
 \end{gather} Since
$x^{(n)}z^{n+1}$ is $O(z)$, the singular part of (\ref{2.27})
equals the right hand side of (\ref{2.26a}).

We now rewrite the right hand side of (\ref{2.26a}) as
\begin{gather*}
\left(\big(x^{(n)}+\zeta_{n+1}z^{-n-1}\big)\sum_{p=0}^n\big(\bar\zeta_{n+1}
x^{(n)}z^{n+1}\big)^p\right)_{-}\\
\qquad{}= \left(\sum_{p=0}^n\bar\zeta_{n+1}^p\big(x^{(n)}\big)^{p+1}z^{p(n+1)}+\zeta_{
n+1}z^{-n-1}+\sum_{p=1}^n\vert\zeta_{n+1}\vert^2\bar\zeta_{n+1}^{
p-1}\big(x^{(n)}\big)^pz^{(p-1)(n+1)}\right)_{-}\\
\qquad\quad{}+\left(\sum_{p=0}^n\bar\zeta_{n+1}^p\big(x^{(n)}\big)^{p+1}z^{p(n+1)}+\zeta_{
n+1}z^{-n-1}+\sum_{p=1}^n\vert\zeta_{n+1}\vert^2\bar\zeta_{n+1}^p
\big(x^{(n)}\big)^{p+1}z^{p(n+1)}\right)_{-}\\
\qquad\quad{}+\left(\sum_{p=0}^n\bar\zeta_{n+1}^p\big(x^{(n)}\big)^{p+1}(1+\vert\zeta_{
n+1}\vert^2)z^{p(n+1)}\right)_{-}+\zeta_{n+1}z^{-n-1}.
\end{gather*} This
completes the proof of part (a).

Part (b) is obvious.\end{proof}

For small $n$ the recursion implies $x^{(1)}=\zeta_1z^{-1}$,
\begin{gather}\label{x2}x^{(2)}=\zeta_1\big(1+\vert\zeta_2\vert^2\big)z^{-1}+\zeta_2z^{-2},\\
\label{x3}
x^{(3)}=\big(\zeta_1\big(1+\vert\zeta_2\vert^2\big)\big(1+\vert\zeta_3\vert^2\big)+
\zeta_2\big(1+\vert\zeta_3\vert^2\big)\zeta_2\bar{\zeta}_3\big)z^{-1}
+\zeta_2\big(1+\vert\zeta_3\vert^2\big)z^{-2}+\zeta_3z^{-3},\\
x^{(4)}=\left(\zeta_1\prod_2^4\big(1+\vert\zeta_j\vert^2\big)+\zeta_2\prod_3^
4\big(1+\vert\zeta_j\vert^2\big)\big(\zeta_2\bar{\zeta}_3+2\zeta_3\bar{\zeta}_
4\big)+\zeta_3\big(1+\vert\zeta_4\vert^2\big)\big(\zeta_3\bar{\zeta}_4\big)^2\right)z^{-1}\nonumber\\
\phantom{x^{(4)}=}{} +\left(\zeta_2\prod_3^4\big(1+\vert\zeta_j\vert^2\big)+\zeta_3\big(1+\vert\zeta_
4\vert^2\big)\zeta_3\bar{\zeta}_4\right)z^{-2}+\zeta_3\big(1+\vert\zeta_4\vert^2\big)z^{-3}+\zeta_4z^{-4}.
\label{x4}
\end{gather}

\begin{lemma}\label{2.40} $y^{(n)}=(zx^{(n+1)})_{-}$ depends only on
$\zeta_2,\dots,\zeta_{n+1}$, and satisfies the same recursion as~$x^{
(n)}$, with the shifted variables $\zeta_2,\dots$ in place of
$\zeta_1,\dots$.
\end{lemma}

\begin{proof} For small $n$ the formulas above show that
$y^{(n)}$ does not depend on $\zeta_1$.  By Lemma \ref{2.25}
\begin{gather*}
y^{(n)}=\left(z\left(\big(x^{(n)}+\zeta_{n+1}z^{-n-1}\big)\sum_{p=0}^n\big(\bar\zeta_{
n+1}x^{(n)}z^{n+1}\big)^p\right)_{-}\right)_{-}\\
\phantom{y^{(n)}}{} =\left(\big(zx^{(n)}+\zeta_{n+1}z^{-n}\big)\sum_{p=0}^n\big(\bar\zeta_{n+1}
zx^{(n)}z^n\big)^p\right)_{-}\\
\phantom{y^{(n)}}{}
=\left(\big(y^{(n-1)}+\zeta_{n+1}z^{-n}\big)\sum_{p=0}^n\big(\bar\zeta_{n+1}
y^{(n-1)}z^n\big)^p\right)_{-}.
\end{gather*} This establishes the recursion and
induction implies $y^{(n)}$ does not depend on $\zeta_1$.
\end{proof}

We can now complete the proof of part (b). Lemma \ref{2.40}
implies that
\[\big(zx^{(n+1)}\big)_{-}=\sum_2^{n+1}x^{(n+1)}_j(\zeta_1,\dots,\zeta_{n+1}
)z^{-j+1}=\sum_1^nx^{(n)}_i(\zeta_2,\dots,\zeta_n)z^{-i}.
\] This
implies that for $j>1$,
\[
x^{(n+1)}_j(\zeta_1,\dots,\zeta_{n+1})=x^{(n)}_{j-1}(\zeta_2,\dots,\zeta_{
n+1}).\] By induction this implies part (b). This also implies
\[x^{(n+1)}=x^{(n+1)}_1z^{-1}+x^{(n)}(\zeta_2,\dots,\zeta_{n+1})z^{-
1}.\]

For future reference, note that there is a recursion for $x_1$ of
the form
\begin{gather}
x_1(\zeta_1,\dots,\zeta_{n+1})=x_1(\zeta_1,\dots,\zeta_n)\big(1+\vert\zeta_{n+1}\vert^2
\big)\nonumber\\
\qquad {}+\left(\sum_{i+j=n+2}x_1(\zeta_i,\dots,\zeta_n)x_
1(\zeta_j,\dots,\zeta_n)\right)\bar{\zeta}_{n+1}\big(1+\vert\zeta_{n+1}\vert^2
\big)\label{x1recursion}
\\
\qquad {}+\left(\sum_{i+j+k=2n+3}x_1(\zeta_i,\dots,\zeta_n)x_1(\zeta_j,\dots,\zeta_
n)x_1(\zeta_k,\dots,\zeta_n)\right)\bar{\zeta}_{n+1}^2\big(1+\vert\zeta_{n+1}\vert^
2\big)+\cdots.\nonumber
\end{gather}
 It would be highly desirably to f\/ind a closed form
solution of this recursion for $x_1$.

We now consider part (c).  For $g$ as in part (c), consider the
Riemann--Hilbert factorization $g=g_{-}g_0g_{+}$, where
\[g_{-}=\left(\begin{matrix} 1&x\\
0&1\end{matrix} \right),\qquad g_0=\left(\begin{matrix} a_0&b_0\\
0&a_0^{-1}\end{matrix} \right)\in SL(2,\mathbb C ),\] and
$g_{+}\in H^0(D,0;SL(2,\mathbb C ),1)$.  Then
\begin{gather}\label{Zequation}Z(g)=C(g)A(g)^{-1}=C(g_{-})A(g_0g_{+})A(g_0g_{+})^{-1}A(g_{-})=
Z(g_{-}).\end{gather} (For use in the next paragraph, note that
this calculation does not depend on the specif\/ic form of $g_0$.)
Let $\epsilon_1$, $\epsilon_2$ denote the standard basis for
$\mathbb C^2$. As in \cite{PS}, consider the ordered basis
$\dots,\epsilon_1z^{j+1},\epsilon_2z^{j+1},\epsilon_1z^j,\dots$,
$j\in\mathcal Z$, for $\mathcal H$. This basis is compatible with
the Hardy polarization of $\mathcal H$.  We claim that
\begin{gather}\label{6.73}Z(g_{-})=C(g_{-})=\left(\begin{matrix} .&0&x_n&.&0&x_3&0&x_2&0&x_1\\
.&0&0&.&0&0&0&0&0&0\\
.&.&.&&&x_4&0&x_3&0&x_2\\
&&&&&0&0&0&0&0\\
.&&&&&&&&0&x_3\\
.&&.&&&.&.\\
.&&&&&.&&&0&x_n\\
0&0&0&0&.&.&&0&0&0\end{matrix} \right). \end{gather}

Let $P_{\pm}$ denote the orthogonal projections associated to the
Hardy splitting of $\mathcal H$. For example
\begin{gather}\label{hardypolarization}P_{+}\left(f=\sum f_kz^k\right)=\sum_{k\ge 0}f_kz^k. \end{gather}
Suppose that $\left(\begin{matrix} f_1\\
f_2\end{matrix} \right)\in \mathcal H^{+}$.  Then
\begin{gather*}
C(g_{-})A(g_{-})^{-1}\left(\begin{matrix} f_1\\
f_2\end{matrix} \right)=P_{-}g_{-}P_{+}g_{-}^{-1}\left(\begin{matrix} f_1\\
f_2\end{matrix} \right)=P_{-}g_{-}\left(\begin{matrix} f_1-P_{+}(xf_2)\\
f_2\end{matrix} \right),\nonumber\\
P_{-}\left(\begin{matrix} f_1+P_{+}(xf_2)+xf_2\\
f_2\end{matrix} \right)=
\left(\begin{matrix} P_{-}(xf_2)\\
0\end{matrix} \right)=C(g_{-})\left(\begin{matrix} f_1\\
f_2\end{matrix} \right).
\end{gather*} This is the f\/irst part of
the claim. For the second part one simply calculates directly,
using the simple form for $g_{-}$.

Comparing (\ref{6.73}) with (\ref{hank1}) proves the f\/irst part of~(c).

Using the factorization
\[g=\left(\begin{matrix} 1&x\\
0&1\end{matrix} \right)g_0g_{+},\] and the specif\/ic form of $g_0$,
it is clear that
\[g_{-}\left(\left(\begin{matrix} z^{1/2}&0\\
0&z^{-1/2}\end{matrix} \right)g\left(\begin{matrix} z^{-1/2}&0\\
0&z^{1/2}\end{matrix} \right)\right)=\left(\begin{matrix} 1&x'\\
0&1\end{matrix} \right),\] where
$x'=(zx)_{-}=x_2z^{-1}+\cdots+x_nz^{-(n-1)}$.  We now use
(\ref{sigma1}) and (\ref{6.73}) to prove the second part of (c).

Part (d) follows from (c). Part (e) can be read of\/f from Lemma
\ref{2.25} (see (\ref{x2})--(\ref{x4})), or from part~(2) of
Proposition~\ref{Luprop}.
\end{proof}

We are now in a position to prove Theorem \ref{hankeltheorem} at
the end of the Introduction.

\begin{proof} By parts (d)
and (e) of Theorem \ref{LSU(2)case},
\begin{gather*}
\int\frac 1{\prod\limits_{l=0}^{n-1}\det\left(1+B\left(\sum\limits_{j=1}^{n-l}x_{l+j}z^j
\right)B_n\left(\sum\limits_{j=1}^{n-l}x_{l+j}z^j\right)^{*}\right)^{p_l}}d\lambda
(x_1,\dots,x_n)\\
\qquad {}=\int \left(\prod_{l=0}^{n-1}\prod_{j=1}^{n-l}\big(1+\vert\zeta_{l+j}\vert^
2\big)^{-jp_l}\right)\prod_{j=1}^n\big(1+\vert\zeta_j\vert^2\big)^{-2(j-1)}d\lambda
(\zeta_j)\\
\qquad{}=\int \big(1+\vert\zeta_1\vert^2\big)^{-p_1}d\lambda (\zeta_1)\int \big(1+\vert
\zeta_2\vert^2\big)^{2-(2p_1+p_2)}d\lambda (\zeta_2)\cdots\\
\qquad\quad{}\cdots \int \big(1+\vert\zeta_n\vert^2\big)^{2n-2-(np_1+\cdots+p_n)}d\lambda (\zeta_
n)\\
\qquad{} =\pi^n\frac 1{p_1-1}\frac 1{2p_1+p_2-3}\cdots \frac 1{np_1+(n-1)p_2+
\cdots +p_n-(2n-1)}.\tag*{\qed}
\end{gather*}\renewcommand{\qed}{}
\end{proof}

\section{The SU(2) case. II}\label{SU(2)caseII}

This is a continuation of Section \ref{SU(2)case}. We f\/irst
consider the limit $n \to \infty $, in the context of Theorem
\ref{LSU(2)case}. From the point of view of analysis, this limit
is naturally related to the critical exponent $s=1/2$ for the
circle. We secondly show that there is a global factorization of
the momentum mapping for $(S(1),\omega)$, extending the formulas
in (c) of Theorem \ref{LSU(2)case}. As in Section \ref{SU(2)case},
we will continue to view $S(1)$ as a submanifold of $LSU(2)$,
rather than $\tilde{L}SU(2)$.

As in (\ref{hardypolarization}), we let $P_{\pm}$ denote the
orthogonal projections associated to the Hardy splitting of
$L^2(S^1)$. Given $f=\sum f_nz^n$, we will write $f^{*}=\sum\bar
{c}_nz^{-n}$. If we simply write $z$ for the multiplication
operator corresponding to $z$, then
\[(\cdot )^{*}\circ P_{-}=z\circ P_{+}\circ z^{-1}\circ (\cdot )^{
*},\] and
\[(\cdot )^{*}\circ P_{+}=z\circ P_{-}\circ z^{-1}\circ (\cdot )^{
*}.\] For a function $F \in L^{\infty}(S^1)$, viewed as a bounded
multiplication operator on $L^2(S^1)$, we will write
$A(F)=P_+FP_+$, and so on, as in (\ref{segalmatrixnotation}).

Suppose
\[l=\left(\begin{matrix} 1&x\\
0&1\end{matrix} \right),\qquad {\rm where}\quad x=\sum_1^nx_jz^{-j}.\]
There exists a unique $g\in L_{\rm pol}SU(2)$ with unique triangular
factorization $g=la(g)^{h_1}u$ where $a(g)=\vert\sigma_1\vert
(g)/\vert\sigma_0\vert (g)$, and
\begin{gather}\label{detformula}\vert\sigma_0\vert^2=\det(1+C(x)^{*}C(x))^{-1},\qquad\vert\sigma_
1\vert^2=\det(1+C(zx)^{*}C(zx))^{-1}\end{gather} (see part (c) of
Theorem~\ref{LSU(2)case}).

\begin{lemma}\label{factorizationlemma}The triangular factorization of $g\in
L_{\rm pol}SU(2)$ is given by
\[g=\left(\begin{matrix} 1&x\\
0&1\end{matrix} \right)\left(\begin{matrix} a&0\\
0&a^{-1}\end{matrix} \right)\left(\begin{matrix}
1+\sum\limits_1^{n-1}\alpha_jz^j&
\sum\limits_0^{n-2}\beta_jz^j\\
\sum\limits_1^n\gamma_jz^j&1+\sum\limits_1^{n-1}\delta_jz^j\end{matrix}
\right),
\] where
\begin{gather*}
\gamma =-((1+C(zx)C(zx)^{*})^{-1}(x))^{*}, \qquad \delta^*=C(x)\gamma,\\
 1+\alpha =\frac 1{a^2}(1-A(x)\gamma) , \qquad \beta =-\frac 1{a^2}A(x)(1+\delta ).
\end{gather*}
\end{lemma}

\begin{proof} Because
\[g=\left(\begin{matrix} a(1+\alpha )+a^{-1}x\gamma&a\beta +xa^{-1}(1+\delta
)\\
a^{-1}\gamma&a^{-1}(1+\delta )\end{matrix} \right)\] has values in
$SU(2)$ (as a function of $z \in S^1$),
\[a^2(1+\alpha )+x\gamma =1+\delta^{*},\qquad \mbox{and}\qquad a^2\beta +x
(1+\delta )=-\gamma^{*}.\] The f\/irst equation can be expressed in
operator language as
\begin{gather}\label{eqn1}A(x)\gamma =1-a^2(1+\alpha ),\qquad
C(x)\gamma =\delta^{*}. \end{gather} The second equation is
equivalent to
\begin{gather}\label{eqn2}A(x)(1+\delta )=-a^2\beta ,\qquad C(x)(\delta
)=-x-\gamma^{*}.\end{gather} We can solve for $\delta$, using
the second equation of (\ref{eqn1}),
\[(\cdot )^{*}\circ C(x)(\gamma )=zB(z^{-1}x^{*})(\gamma^{*})=\delta\]
The second equation of (\ref{eqn2}) now implies
\[C(x)zB\big(z^{-1}x^{*}\big)(\gamma^{*})=-x-\gamma^{*}\]
which is equivalent to
\[(1+C(zx)C(zx)^{*})(\gamma^{*})=-x.\]
The Lemma follows from these equations.
\end{proof}

Let $W^{1/2}$ denote the Sobolev space of (Lebesgue equivalence
classes of) functions having half of a derivative, i.e.\ if $f=\sum
f_j z^j$, then $\sum j \vert f_j \vert ^2 <\infty $. A class in
$W^{1/2}$ is not in general represented by a continuous function.
Despite this, $W^{1/2}(S^1,SU(2))$ is a connected topological
group (homotopy equivalent to $LSU(2)$; see \cite{Brezis}), and it
is the natural domain for the basic factorization theorems in the
theory of loop groups (see Chapter~8 of~\cite{PS}). We let
$\mathbb S(1)$ denote the completion of $S(1)$ in
$W^{1/2}(S^1,SU(2))$.

\begin{theorem}\label{limittheorem} By taking a limit as $n \to \infty$ in Theorem~{\rm \ref{LSU(2)case}},
we obtain bijective correspondences among the following three sets

{\rm (a)} $\big\{\zeta
=(\zeta_1,\zeta_2,\dots):\sum\limits_1^{\infty}j\vert\zeta_j\vert^ 2<\infty
\big\}$.

{\rm (b)} ${x=\sum\limits_1^{\infty}x_jz^{-j}\in W^{1/2}}$, where (as in {\rm (b)} of
Theorem~{\rm \ref{LSU(2)case}})
\begin{gather}\label{xjx1relation}x_j(\zeta_1,\dots)=\lim_{n\to\infty}x_1(\zeta_j,\dots,\zeta_n),\end{gather}
and conversely
\begin{gather}\label{zetarelations}\zeta_j(x_1,\dots)=\lim_{n\to\infty}\zeta_1(x_j,\dots,x_n).\end{gather}

{\rm (c)}
\[\left\{g=\left(\begin{matrix} d(z)^{*}&-c(z)^{*}\\
c(z)&d(z)\end{matrix} \right)\in \mathbb S(1):c(z) ,d(z) \in
H^0(\Delta ),c(0)=0\right\},\] where
\begin{gather*}
g=\lim_{n\to\infty}a(\zeta_n)\left(\begin{matrix} 1&\zeta_nz^{-n}\\
-\bar{\zeta}_nz^n&1\end{matrix}
\right)\cdots a(\zeta_1)\left(\begin{matrix}
1&\zeta_1z^{-1}\\
-\bar{\zeta}_1z&1\end{matrix} \right),\end{gather*}
$a(\zeta_j)=(1+\vert \zeta_j\vert^2)^{-1/2}$, $g$ has triangular
factorization
\begin{gather}\label{trifac}g=\left(\begin{matrix} 1&\sum\limits_1^{\infty}x_jz^{-j}\\
0&1\end{matrix} \right)\left(\begin{matrix} \frac
{\vert\sigma_1\vert (g)}{
\vert\sigma_0\vert (g)}&0\\
0&\frac {\vert\sigma_0\vert (g)}{\vert\sigma_1\vert
(g)}\end{matrix} \right )u,\end{gather} and the entries of $u$
are given by the same formulas as in Lemma~{\rm \ref{factorizationlemma}}.
\end{theorem}

Note that for $k>0$ $g=\left(
\begin{matrix}z^{-k}&0\\0&z^k\end{matrix} \right)$ has a second
row which is holomorphic in the disk, but $g$ is not in $S(1)$.
Thus in part (c) it is necessary to require that $g \in \mathbb
S(1)$.

\begin{proof}  By Theorem \ref{LSU(2)case}, for $\zeta$ with a f\/inite number of terms,
\begin{gather}\label{product}\prod \big(1+\vert\zeta_j\vert^2\big)^j=\det(1+C(x)C(x)^{*}),\end{gather}
and obviously
\begin{gather}\label{trace}{\rm tr}(C(x)C(x)^{*})=\sum j\vert x_j\vert^2.\end{gather}
For an arbitrary sequence $\zeta$, the product on the LHS of
(\ref{product}) is f\/inite if\/f $\sum\limits_1^{\infty}j \vert \zeta_j
\vert^2<\infty$. Similarly, for an arbitrary $x$, the determinant
is f\/inite if\/f (\ref{trace}) is f\/inite. Thus given $\zeta$ as in
part~(a), the partial sums for the series representing $x$ will
have limits in $W^{1/2}$. To understand why there is a unique
limit, and to prove the other statements in part (b), recall that
the recursion relation~(\ref{x1recursion}) implies that $x_1$ has
a series expansion, with nonnegative integer coef\/f\/icients, in
terms of the variables $\zeta_j$, $\bar{\zeta_j}$, of the form
\begin{gather*}
x_1(\zeta_1,\dots)=\zeta_1 \prod_{j=2}^{\infty}(1+\vert\zeta_j\vert^2)
+\zeta_2\prod_{j=3}^{\infty}\big(1+\vert\zeta_j\vert^2\big)\big(\zeta_2\bar{\zeta}_3+2\zeta_3\bar{\zeta}_
4+\cdots\big)\nonumber\\
\phantom{x_1(\zeta_1,\dots)=}{} +\zeta_3\prod_{j=4}^{\infty}\big(1+\vert\zeta_j\vert^2\big)\big((\zeta_3\bar{\zeta}_4)^2+\cdots\big)+\cdots,
\end{gather*}
Since $\vert x_1 \vert^2$ is dominated by (\ref{trace}), this
series will converge absolutely. Thus $x_1$ is a well-def\/ined
function of $\zeta$. The relation (\ref{xjx1relation}) follows
from (b) of Theorem (\ref{LSU(2)case}). Thus all the $x_j$, and
hence also $x \in W^{1/2}$, are uniquely determined by $\zeta$,
assuming $\sum j\vert \zeta_j \vert^2<\infty$.

Because
\[\prod_{j=k}^{\infty}\big(1+\vert\zeta_j\vert^2\big)=\frac{\det\big(1+C(z^kx)C(z^kx)^*\big)}
{\det\big(1+C(z^{k-1}x)C(z^{k-1}x)^*\big)}\] and the triangular nature of
the relation between the $\zeta_j$ and the $x_j$, the map from the
$\zeta$ to $x$ can be inverted, and the $\zeta_j$ will be
expressible in terms of $\zeta_1$ as in (\ref{zetarelations}). We
have thus established that there is a bijective correspondence
between the sets in (a) and (b).

Now suppose that we are given $x$ as in part (b). We claim that
(\ref{detformula}) and Lemma \ref{factorizationlemma} imply that
we can obtain a $g$ as in part (c), with triangular factorization
as in (\ref{trifac}). Because $x\in W^{1/2}$, $C(x)$ and $C(zx)$
are Hilbert--Schmidt operators (viewed as operators on $L^2(S^1)$),
so that the formulas~(\ref{detformula}), and hence the formula for
$a(g)$, make sense. The formulas for $\gamma$ and $\delta$ in
Lemma~\ref{factorizationlemma} a priori show only that $\gamma$
and $\delta$ are $L^2(S^1)$ (not necessarily $W^{1/2}$), so that
we obtain a $SU(2)$ loop $g$, expressed as in part (c), with $L^2$
entries. But as in (\ref{Zequation}),
$Z(g)=Z\left(\left(\begin{matrix}1&x\\0&1\end{matrix}\right)\right)$, and
this is a Hilbert--Schmidt operator. Because $g \in
L^{\infty}(S^1)$, $A(g)$ is a bounded operator, and hence this
implies that $C(g)$ is a Hilbert--Schmidt operator. Thus $g \in
W^{1/2}$.

Conversely given $g$ as in part (c), it is obvious that $g$ has a
triangular form as in (\ref{trifac}), and this determines $x$. The
equality
\[\det\vert A(g)\vert^2=\det\big(1+C(x)C(x)^{*}\big)^{-1}\]
implies that $x\in W^{1/2}$. Because the sets in (a) and (b) are
in correspondence, this also determi\-nes~$\zeta$. This completes
the proof.
\end{proof}

Recall that in Theorem \ref{LSU(2)case} we considered the minimal
sequence $\alpha_0,\alpha_1,\dots$. By considering the sequence
$\alpha_1,\alpha_0,\dots$, or in other words, by applying the
automorphism (\ref{automorphism}) which interchanges the two
simple roots, the proceeding Theorem can be reformulated in the
following way.

\begin{corollary}\label{limittheorem2} There are bijective correspondences among the following three sets

{\rm (a)} $\big\{\eta =(\eta_0,\eta_1,\dots):\sum\limits_0^{\infty}j\vert\eta_j\vert^
2<\infty \big\}$.

{\rm (b)} ${y=\sum\limits_0^{\infty}y_jz^{-j}\in W^{1/2}}$, where
\begin{gather*}
y_j(\eta_0,\dots)=\lim_{n\to\infty}y_0(\eta_j,\dots,\eta_n),
\end{gather*}
and conversely
\begin{gather*}
\eta_j(y_0,\dots)=\lim_{n\to\infty}\eta_0(y_j,\dots,y_n).\end{gather*}

{\rm (c)}
\[\left\{h=\left(\begin{matrix} e(z)&b(z)\\
-b^*(z)&e(z)^*\end{matrix} \right)\in \mathbb S(1):b(z),e(z) \in
H^0(\Delta )\right\},\] where
\begin{gather*}
h=\lim_{n\to\infty}a(\eta_n)\left(\begin{matrix} 1&\eta_nz^n\\
-\bar{\eta}_nz^{-n}&1\end{matrix}
\right)\cdots a(\eta_0)\left(\begin{matrix}
1&\eta_0\\
-\bar{\eta}_0&1\end{matrix} \right),\end{gather*} and $h$ has
triangular factorization of the form
\begin{gather}\label{htrifac}h=\left(\begin{matrix} 1&0\\
y&1\end{matrix} \right)\left(\begin{matrix} \frac
{\vert\sigma_1\vert (h)}{
\vert\sigma_0\vert (h)}&0\\
0&\frac {\vert\sigma_0\vert (h)}{\vert\sigma_1\vert
(h)}\end{matrix} \right ) \left(\begin{matrix}
1+\sum\limits_1^{\infty}\alpha_jz^j&
\sum\limits_0^{\infty}\beta_jz^j\\
\sum\limits_1^{\infty}\gamma_jz^j&1+\sum\limits_1^{\infty}\delta_jz^j\end{matrix}
\right),\\
\beta =-(1+C(y)^*C(y))^{-1}(y^*), \qquad \alpha^*=C(y)\beta,\nonumber\\
 1+\delta =a(h)^2(1-A(y)\beta) , \qquad \gamma^*=-a(h)^2D(y^*)(\alpha^*).\nonumber
\end{gather}
Moreover
\[\vert\sigma_0\vert^2(h)=\prod_{j=0}^{\infty}\big(1+\vert\eta_j\vert^
2\big)^{-j},\qquad \vert\sigma_1\vert^2(h)=\prod_{j=0}^{\infty}\big(1+\vert\eta_j\vert^
2\big)^{-(j+1)}\] and
\[a(h)=\frac {\vert\sigma_1\vert}{\vert\sigma_0\vert}=\prod \big(1+\vert
\eta_j\vert^2\big)^{-1/2}.\]
\end{corollary}

In the following statement we will continue to view $\mathbb S(1)$
as a subset of $W^{1/2}(S^1,SU(2))$, but in the proof it will be
necessary to consider lifts in the Kac--Moody extension (as in
Lemma \ref{notationlemma}).

\begin{theorem} Suppose that $\zeta$, $\chi$, and $\eta$ are sequences such that
\[
\sum_{j=1}^{\infty}j\big(\vert\zeta_j\vert^2+\vert\chi_j\vert^2+\vert\eta_j\vert^2\big)<\infty.\]
By slight abuse of notation, we identify $\chi$ with the function
$\chi=\sum\limits_{j=1}^{\infty}\chi_jz^j$. Let $g$ be defined as in~{\rm (c)}
of Theorem~{\rm \ref{limittheorem}} and $h$ as in {\rm (c)} of Corollary~{\rm \ref{limittheorem2}}.

{\rm (a)} The product $h^{-1}e^{(\chi -\chi^{*})h_1}g\in \mathbb S(1)$.

{\rm (b)} The mapping
\[
{(\eta,\chi,\zeta): \ \sum j \big (\vert \eta_j \vert^2+ \vert \chi
_j \vert^2+ \vert \zeta_j \vert^2\big )<\infty} \to \mathbb
S(1): \ (\eta,\chi,\zeta)\to h^{-1}e^{(\chi -\chi^{*})h_1}g \] is
injective.

{\rm (c)}
\begin{gather*}
\vert\sigma_0\vert^2\big(h^{-1}e^{(\chi -\chi^{*})h_1}g\big)=\prod_{j=1}\big(1+\vert
\eta_j\vert^2\big)^{-j}\exp\left(-2\sum_{j=1}j\vert\chi_j\vert^2\right)\prod_{j=1}
\big(1+\vert\zeta_j\vert^2\big)^{-j},\\
\vert\sigma_1\vert^2\big(h^{-1}e^{(\chi -\chi^{*})h_1}g\big)=\prod_{j=1}(1+\vert
\eta_{j+1}\vert^2)^{-j}\exp\left(-2\sum_{j=1}j\vert\chi_j\vert^2\right)\prod_{
j=1}\big(1+\vert\zeta_{j+1}\vert^2\big)^{-j},\\
a^2\big(h^{-1}e^{(\chi -\chi^{*})h_1}g\big)=\prod_{j=1}\frac {1+\vert\zeta_j
\vert^2}{1+\vert\eta_j\vert^2}.
\end{gather*}
\end{theorem}

\begin{proof}For part (a), we do a calculation at the level of loops.
If we write $h$ as in (\ref{htrifac}), then the triangular
factorization of $h^{-1}$ is given by
\begin{gather*}
h^{-1}=l\big(h^{-1}\big)\left(\begin{matrix} a(h^{-1})&0\\
0&a(h^{-1})^{-1}\end{matrix} \right)\left(\begin{matrix} 1&y^*\\
0&1\end{matrix} \right),\end{gather*} where
\begin{gather*}
l\big(h^{-1}\big)=u(h)^*=\left(\begin{matrix} 1+\alpha^{*}&\gamma^{*}\\
\beta^{*}&1+\delta^{*}\end{matrix}\right).
\end{gather*} Then
\begin{gather}\label{firsteqn}h^{-1}e^{(\chi -\chi^{*})h_1}g=l\big(h^{-1}\big)a\big(h^{-1}\big)^{h_1}\left(\begin{matrix} 1&y^*\\
0&1\end{matrix} \right)e^{(\chi -\chi^{*})h_1}\left(\begin{matrix} 1&x\\
0&1\end{matrix} \right)a(g)^{h_1}u(g).\end{gather} The main
point of the proof is that the middle three factors are upper
triangular, and we can f\/ind the (loop space) triangular
factorization of the product with ease. Thus (\ref{firsteqn})
equals
\begin{gather*}
l\big(h^{-1}\big)e^{-\chi^{*}h_1}\left(\begin{matrix} 1&a(h^{-1})^2P_-(y^*e^{2\chi^{*}}+x
e^{2\chi})\\
0&1\end{matrix} \right)\\
\qquad{}\times a(h)^{h_1}a(g)^{h_1}\left(\begin{matrix}
1&a(g)^{-2}P_+(y^*e^{2\chi^{*}}+xe^{2\chi})\\
0&1\end{matrix} \right)e^{\chi h_1}u(g).
\end{gather*} This triangular form
implies part (a) of the Theorem.

This calculation also implies that $a(h^{-1}e^{(\chi
-\chi^{*})h_1}g)=a(h^{-1})a(g)$. Given that
$a(h^{-1})=a(h^*)=a(h)$, and formulas we have already established
for $a(g)$ and $a(h)$ (see Corollary \ref{limittheorem2}), this
also implies the third formula in part (c).

To obtain the f\/irst two formulas in part (c), we need to lift each
of the factors, $h^{-1}$, the torus-valued loop, and $g$, into the
extension (where the lift is determined by requiring that the lift
is in $\mathbb S(1)$), and then repeat the preceding calculation
in the extension. To do this we replace~$a(g)^{h_1}$ by
$\vert\sigma_0\vert(g)^{h_0}\vert\sigma_1\vert(g)^{h_1}$, and
similarly for $h^{-1}$ (and recall that we have explicit product
formulas for the functions $\vert\sigma_j\vert(g)$, and similarly
for $h$). The torus-valued loop $e^{(\chi-\chi^*)h_1}$ has a
vanishing diagonal term in its triangular factorization. However,
a well-known formula for Toeplitz determinants implies that
\[\vert\sigma_j\vert\big(e^{(\chi-\chi^*)h_1}\big)=\exp\left(-2\sum_{j=1}j\vert\chi_j\vert^2\right).\]
Thus the lift of this torus-valued loop into the extension has
diagonal term
\[\exp\left(-2\sum_{j=1}j\vert\chi_j\vert^2(h_0+h_1)\right).\]
It is now straightforward to repeat the calculation above and
conclude, as we did at the level of loops, that (i) the product of
these lifts is in $\mathbb S(1)$ and (ii) the diagonal term of the
product is the product of the diagonal terms of the factors. This
implies the formulas in part~(c).

To prove that the mapping is injective, recall that
\begin{gather}
l\big(h^{-1}e^{(\chi -\chi^{*})h_1}g\big)=l\big(h^{-1}\big)e^{-\chi^{*}h_1}a(h)^{h_1}\left
(\begin{matrix} 1&P_-(y^*e^{2\chi^{*}}+xe^{2\chi})\\
0&1\end{matrix} \right)a(h)^{-h_1}\nonumber\\
\qquad {}=\left(\begin{matrix} 1+\alpha^{*}&\gamma^{*}\\
\beta^{*}&1+\delta^{*}\end{matrix}
\right)e^{-\chi^{*}h_1}\left(\begin{matrix}
1&a(h)^2P_-(y^*e^{2\chi^{*}}+xe^{2\chi})\vspace{1mm}\\
0&1\end{matrix} \right)\nonumber\\
\label{lmatrix}
\qquad{}=\left(\begin{matrix}
(1+\alpha^{*})e^{-\chi^{*}}&(1+\alpha^{*})e^{-\chi^{*}}a
(h)^2P_-(y^*e^{2\chi^{*}}+xe^{2\chi})+\gamma^{*}e^{\chi^{*}}\vspace{1mm}\\
\beta^{*}e^{-\chi^{*}}&\beta^{*}e^{-\chi^{*}}a(h)^2P_-(y^*e^{2\chi^{
*}}+xe^{2\chi})+(1+\delta^{*})e^{\chi^{*}}\end{matrix}
\right).\end{gather} We need to show that this matrix determines
$y$ (equivalently $h$ or $\eta$), $\chi$, and $x$ (or $g$ or
$\zeta$).

By the form of the triangular factorization of $h$ in part (c) of
Corollary~\ref{limittheorem2}, it is clear that the f\/irst row of
$h$, evaluated at a point $z\in S^1$, is determined by
$1+\alpha(z)$ and $\beta(z)$. Because $h(z)\in SU(2)$, this also
determines $h(z)$. This means that the f\/irst row of $h$, as a
holomorphic function of $z\in \Delta$ is determined up to a phase
ambiguity by the ratio $\beta/(1+\alpha)$.

Now suppose that we are given (\ref{lmatrix}). The f\/irst column
determines the pair of holomorphic functions $\beta$ and
$1+\alpha$ up to multiplication by a holomorphic function. This
determines the ratio $\beta/(1+\alpha)$, and it also f\/ixes the
phase ambiguity. Thus the f\/irst column determines $h$. It is then
clear that $\chi$ is determined. We can then use the $(1,2)$ entry
to f\/ind $x$. This proves injectivity, and completes the proof of
the Theorem.
\end{proof}

\appendix

\section{Appendix}\label{AppendixA}

In this appendix we will review some of the ideas in \cite{Lu},
relevant to this paper, from a slightly dif\/ferent perspective. The
main rationale for including this appendix is that the basic
arguments are valid in the more general Kac--Moody category.
Throughout this appendix, we will use the notation and basic
results in \cite{Kac}.

We start with the following data:  $A$ is an irreducible
symmetrizable generalized Cartan matrix; $\mathfrak g=\mathfrak
g(A)$ is the corresponding Kac--Moody Lie algebra, realized via its
standard (Chevalley--Serre) presentation; $\mathfrak g=\mathfrak
n^{-}\oplus \mathfrak h\oplus \mathfrak n^{+}$ is the triangular
decomposition; $\mathfrak b =\mathfrak h\oplus \mathfrak n^{+}$
the upper Borel subalgebra; $G=G(A)$ is the algebraic group
associated to $A$ by Kac--Peterson; $H$, $N^{\pm}$ and $B$ are the
subgroups of $G$ corresponding to $\mathfrak h$, $\mathfrak
n^{\pm}$, and $\mathfrak b $, respectively; $ K$ is the ``unitary
form'' of~$G$; $T=K\cap H$ the maximal torus; and $W=N_K(T)/T\simeq
N_G(H)/H$ is the Weyl group.

A basic fact is that $(G,B,N_G(H))$ with Weyl group $W$ is an
abstract Tits system.  This yields a complete determination of all
the (parabolic) subgroups between $B$ and $G$.  They are described
as follows.

Let $\Phi$ be a f\/ixed subset of the simple roots.  The subgroup of
$W$ generated by the simple ref\/lections corresponding to roots in
$\Phi$ will be denoted by $W(\Phi )$. The parabolic subgroup
corresponding to $\Phi$, $P=P(\Phi )$, is given by $P=BW(\Phi )B$.
Given $\mathbf w\in N_K(T)$, we will denote its image in $W/W(\Phi
)$ by $w$.

The basic structural features of $G/P$ which we will need are the
Birkhof\/f and Bruhat decompositions
\begin{gather*}
G/P=\bigsqcup\Sigma_w,\qquad\Sigma_w=N^{-}\mathbf w P,
\\
G/P=\bigsqcup C_w,\qquad C_w=B\mathbf w P,
\end{gather*}
respectively, where the indexing set is $W/W(\Phi )$ in both
cases.  The strata $\Sigma_w$ are inf\/inite dimensional if
$\mathfrak g$ is inf\/inite dimensional, while the cells $C_w$ are
always f\/inite dimensional.  Our initial interest is in the
Schubert variety $\bar {C}_w$, the closure of the cell.

Fix $w\in W/W(\Phi )$.  We choose a representative $\mathbf w\in
N_K(T)$ of minimal length $n$, and we f\/ix a factorization
\begin{gather}\label{A3}
\mathbf w=r_n\cdots r_1,
\end{gather}
where $r_j=i_{\gamma_i}\begin{pmatrix} 0&i\\
i&0\end{pmatrix}$, and $i_{\gamma_j}:SL_2\rightarrow G$ is the
canonical homomorphism of $SL_2$ onto the root subgroup
corresponding to the simple root $\gamma_j$.

\begin{proposition}\label{dougA4}
For $\mathbf w$ as in \eqref{A3}, the map
\[
r_n\exp(\mathfrak g_{-\gamma_n})\times \cdots \times r_1\exp(\mathfrak
g_{-\gamma_ 1})\to G/P: \ (p_j)\to p_n \cdots p_1P
\]
is a complex analytic isomorphism onto $C_w$.
\end{proposition}

This result is essentially~(5) of \cite{Kac} together with Tits's
theory.  We will include a proof for completeness.

\begin{proof} Let $\Delta^{+}$ denote the positive roots,
$\Delta^{+}(\Phi )$ the positive roots which are combinations of
elements from $\Phi$.  The ``Lie algebra of $P$'' is $\mathfrak
p=\Sigma \mathfrak g_{-\beta}\oplus \mathfrak b $ where the sum is
over $\beta\in\Delta^{+}(\Phi )$; this is the Lie algebra of $P$
in the sense that it is the subalgebra generated by the root
spaces~$\mathfrak g_{\gamma}$ for which $\exp:\mathfrak g_{\gamma}
\rightarrow G$ is def\/ined and have image contained in $P$. The
subgroups $\exp(\mathfrak g_{\gamma} )$ generate $P$. We also let
$\mathfrak p^{-}$ denote the subalgebra opposite $\mathfrak p$:
$\mathfrak p^{-}=\sum \mathfrak g_{-\gamma}$, where the sum is
over $\gamma\in\Delta^{+}\setminus\Delta^{+}(\Phi )$.  The
corresponding group will be denoted by $P^{-}$.

The cell $C_w$ is the image of the map $N^{+}\rightarrow
G/P:u\rightarrow u\mathbf w P$. The stability subgroup at $\mathbf
w P$ is $N^{+}\cap \mathbf w P\mathbf w^{ -1}$.

At the Lie algebra level we have the splitting
\begin{gather}\label{splitting}
\mathfrak n^{+}=\mathfrak n^{+}\cap {\rm Ad}(\mathbf w)(\mathfrak
p)\oplus \mathfrak n^{+}\cap {\rm Ad}(\mathbf w)(\mathfrak p^{-}).
\end{gather}
The second summand equals $\mathfrak n_w^{+}=\oplus\,\mathfrak
g_{\beta}$, where the sum is over roots $\beta >0$ with
$w^{-1}\beta\in -(\Delta^{+}\setminus\Delta^{+}(\Phi ))$.  These
roots $ \beta$ are necessarily real, so that $\exp:\mathfrak
n_w^{+}\rightarrow N_w^{+}\subseteq N^{+}$ is well-def\/ined.

For $q\in \mathbb Z^{+}$ let $N_q^{+}$ denote the subgroup
corresponding to $\mathfrak n_q^{+}={\rm span}\{\mathfrak
g_{\beta}:{\rm height}(\beta)\geq q\}$.  Then $ N^{+}/N_q^{+}$ is a
f\/inite dimensional nilpotent Lie group, and it is also simply
connected. By taking $q$ suf\/f\/iciently large and considering the
splitting (\ref{splitting}) modulo $\mathfrak n_q^{+}$, we
conclude by f\/inite dimensional considerations that each element in
$N^{+}$ has a unique factorization $n=n_1n_2$, where $n_1\in
N_w^{+}$ and $n_2\in N^{+}\cap \mathbf w P\mathbf w^{-1}$:
\begin{gather}\label{dougA7}
N^{+}\simeq N_w^{+}\times \big(N^{+}\cap \mathbf w P\mathbf w^{-1}\big).
\end{gather}
The important point is that modulo $N_q^{+}$, we can control
$N^{+}\cap \mathbf w P\mathbf w^{-1}$ by the exponential map.

The following lemma is standard.

\begin{lemma}\label{dougA8} In terms of the minimal
factorization $\mathbf w=r_n\cdots r_1$, the roots $\beta >0$ with
$\mathbf w^{-1}\beta <0$ are given by
\[
\beta_j=r_n\cdots r_{j+1}(\gamma_j)=r_n\cdots
r_j(-\gamma_j),\qquad 1\leq j\leq n.
\]
\end{lemma}

Because $\mathbf w$ is a representative of $w\in W/W(\Phi )$ of
minimal length, all of these $\beta_j$ satisfy $\mathbf
w^{-1}\beta_j\in -( \Delta^{+}\setminus\Delta^{+}(\Phi ))$.
Otherwise, if say $\mathbf w^{-1}\beta_j\in -\Delta^{+}(\Phi )$,
then \[ \mathbf w^{-1}r_{\beta_j}\mathbf w=r_1\cdots
r_{j-1}r_jr_{j-1}\cdots r_1\in N_K(T)\cap P \] and $\mathbf
w'=\mathbf w(\mathbf w^{-1}r_{\beta_j}\mathbf w)=r_n\cdots\widehat
{ r}_j\cdots r_1$ would be a representative of $w$ of length $<n$
(here we have used the fact that $W(\Phi )=N_K(T)\cap P/T$, which
follows from the Bruhat decomposition). For future reference we
note this proves that
\begin{gather}\label{dougA10}
N_w^{+}=N^{+}\cap (N^{-})^w=N^{+}\cap (P^{-})^w
\end{gather}
and (\ref{dougA7}) shows that
\begin{gather}\label{dougA11}
N_w^{+}\times w\cong C_w.
\end{gather}
Now for any $1\le p\le q\le n$, $\bigoplus_{p\le j\le p}\mathfrak
g_{\beta_ j}$ is a subalgebra of $\mathfrak n_w^{+}$.  Thus by
(\ref{dougA7}) \[ \exp(\mathfrak g_{\beta_n})\times\cdots\times
\exp(\mathfrak g_{\beta_1}) \times w\cong C_w. \] This completes
the proof of Proposition \ref{dougA4}, when we write
\begin{gather*}
\exp(\mathfrak g_{\beta_j})=r_n\cdots r_j\exp(\mathfrak
g_{\gamma_j})r_j\cdots r_n .\tag*{\qed}
\end{gather*}
  \renewcommand{\qed}{}
\end{proof}

For each $j$, let $P_j$ denote the parabolic subgroup
$i_{\gamma_j}(SL_2)B$.  Let
\[
\gamma_\mathbf w=P_n\times_B\cdots\times_BP_1/B, \] where \[
P_n\times\cdots\times P_1\times B\times\cdots\times B\rightarrow
P_n\times\cdots\times P_1
\]
is given by \[ (p_j)\times (b_j)\rightarrow
\big(p_nb_n,b_n^{-1}p_{n-1}b_{n-1},\dots ,b_2^{-1}p_1b_1\big).
\]
We have written ``$\gamma_\mathbf w$'' instead of ``$\gamma_w$'' to
indicate that this compact complex manifold depends upon the
factorization (\ref{A3}).

\begin{corollary} \label{dougA17}
The map
\[
\gamma_\mathbf w\rightarrow\bar {C}_w: \ (p_j)\rightarrow p_n\cdots
p_1P
\]
is a (Bott--Samelson) desingularization of $\bar {C}_w$.
\end{corollary}

This is an immediate consequence of Proposition \ref{dougA4}.

Let
\[
SL_2'=\left\{g=\begin{pmatrix} a&b\\
c&d\end{pmatrix} \in SL(2,\mathbb C ):a\ne 0\right\}.
\]

\begin{corollary} \label{dougA18}
Consider the surjective map
\[
SL_2\times\cdots\times SL_2\rightarrow\bar {C}_w: \ (g_j)\rightarrow
r_ni_{\gamma_n}(g_n)\cdots r_1i_{\gamma_1}(g_1)P.
\]
The inverse image of $C_w$ is $SL_2'\times\cdots\times SL_2'$.
\end{corollary}

\begin{proof} Let $\sigma =r_{n-1}\cdots r_1$.  It suf\/f\/ices to
show that for the natural actions
\begin{gather}\label{dougA19}
r_ni_{\gamma_n}(SL_2')\times C_{\bar{\sigma}}\rightarrow C_w,
\end{gather}
\[
r_ni_{\gamma_n}(SL_2\setminus SL_2')\times
C_{\bar{\sigma}}\rightarrow \bar {C}_w, \] and
\[
r_ni_{\gamma_n}(SL_2)\times (\bar {C}_{\bar{\sigma}}\setminus C_{
\bar{\sigma}})\rightarrow\bar {C}_w\setminus C_w.
\]
The f\/irst line, (\ref{dougA19}), follows from Proposition~\ref{dougA4}, since $i_{\gamma_n}(SL_2')\!\subseteq\! \exp(-\mathfrak
g_{-\gamma_n})B$ and $B{\times} C_{\bar{\sigma}}\subseteq
C_{\bar{\sigma}}$.  The second line follows from \[
r_ni_{\gamma_n}\left(\begin{matrix} 0&b\\
c&d\end{matrix} \right)\cdot
C_{\bar{\sigma}}=i_{\gamma_n}\left(\begin{matrix}
c&b\\
0&d\end{matrix} \right)\cdot C_{\bar{\sigma}}\subseteq
C_{\bar{\sigma}} .
\]
For the third line it's clear that the image of the left hand side
is a union of cells, since we can replace $r_ni_{\gamma_n}(SL_2)$
by $P_n$.  This image is at most $n-1$ dimensional.  Therefore it
must have null intersection with $C_w$.
\end{proof}

\begin{corollary}\label{dougA18b} \qquad {}

{\rm (a)} Let $k(\zeta)$ be defined as in
\eqref{su2notation}. The map
\[\mathbb C ^n \to C_w:(\zeta_n,\dots,\zeta_1) \to
r_ni_{\gamma_n}(k(\zeta_n))\cdots r_1i_{\gamma_1}(k(\zeta_1))P\]
is a real analytic isomorphism.

{\rm (b)} In terms of the parameterization in {\rm (a)}, and the
parameterization of $C_w$ by $N^{+}\cap \mathbf w P \mathbf w
^{-1}$ (see \eqref{dougA10} and \eqref{dougA11}), Haar measure
(unique up to a constant) is given by
\begin{gather*}
d\lambda_{N^{+}\cap wPw^{-1}}=\prod_{j=1}^n \big(1+\vert\zeta_j\vert^2\big)^{\delta(w_{j-1}^{-
1}h_{\gamma_j}w_{j-1})-1}\\
\phantom{d\lambda_{N^{+}\cap wPw^{-1}}}{} =\prod_{1\le i<j\le
n}\big(1+\vert\zeta_j\vert^2\big)^{-\gamma_i(w_{i-1}
w_{j-1}^{-1}h_{\gamma_j}w_{j-1}w_{i-1}^{-1})},
\end{gather*}
 where $\delta
=\sum\Lambda_j$, the sum of the dominant integral functionals for
$\mathfrak g$.
\end{corollary}

\begin{proof}The proof of (a) is by induction on $n$. We write $w_n$ in place of
$w$, and $k(\zeta_n)=l(\zeta_n)a(\zeta_n)u(\zeta_n)$ for its
$SL(2,\mathbb C )$ triangular factorization.

The case $n=1$ is obvious. Assume the result holds for $n-1$.
Suppose that
\[r_ni_{\gamma_n}(k(\zeta_n))\cdots
r_1i_{\gamma_1}(k(\zeta_1))P=r_ni_{\gamma_n}(k(\zeta_n^{\prime}))\cdots
r_1i_{\gamma_1}(k(\zeta_1^{\prime}))P.\] Since $C_{w_{n-1}}$ is
$B^+$-stable, this equation implies
\[r_ni_{\gamma_n}(l(\zeta_n))xP=r_ni_{\gamma_n}(l(\zeta_n^{\prime}))yP,
\]
where $xP,yP \in C_{w_{n-1}}$. Proposition \ref{dougA4} implies
that $\zeta_n=\zeta_n^{\prime}$. Given this, the induction
hypothesis implies that $\zeta_j=\zeta_j^{\prime}$ for $j<n$. Thus
the map is injective.

Since $k(\zeta_n)=l(\zeta_n)$ modulo $B^+$, and elements of $B^+$
stabilize $C_{w_{n-1}}$, Proposition~\ref{dougA4} also implies the
map is surjective. This proves (a).

Now consider (b). We will f\/irst establish the second formula by
induction. We will then show the two formulas are equivalent.

Recall the factorization
\[i_{\gamma}(k(\zeta))=i_{\gamma}\left(\left(\begin{matrix} 1&0\\
\zeta&1\end{matrix} \right)\left(\begin{matrix} a(\zeta)&0\\
0&a(\zeta)^{-1}\end{matrix} \right)\left(\begin{matrix} 1&-\bar{\zeta}\\
0&1\end{matrix} \right)\right).\]
From this the case $n=1$ is clear. Now
suppose the second formula in (b) is valid for $n-1$ and $n>1$. We
have
\[r_ni_{\gamma_n}(k(\zeta_n))=r_ni_{\gamma_n}\left(\left(\begin{matrix} 1&0\\
\zeta_n&1\end{matrix} \right)\left(\begin{matrix} a(\zeta_n)&0\\
0&a(\zeta_n)^{-1}\end{matrix} \right)\left(\begin{matrix} 1&-\bar{\zeta}_n\\
0&1\end{matrix} \right)\right).\] The action of
\[i_{\gamma_n}\left(\left(\begin{matrix} 1&-\bar{\zeta}_n\\
0&1\end{matrix} \right)\right)\] on $C_{w_{n-1}}$, in terms of the
parameterization by $N^{+}\cap w_{n-1}N^{-}w_{n-1}^{-1}$, is by
translation. So this action preserves Lebesgue measure. The action
of
\[i_{\gamma_n}\left(\left(\begin{matrix} a(\zeta_n)&0\\
0&a(\zeta_n)^{-1}\end{matrix} \right)\right)\] on $C_{w_{n-1}}$, in
terms of the parameterization by $N^{+}\cap
w_{n-1}N^{-}w_{n-1}^{-1}$, is by conjugation, and we can easily
calculate the ef\/fect on volume. In a routine way this leads to the
second formula in~(b).

To prove the formulas in (b) are equivalent, recall that for a
simple positive root $\gamma$, with corresponding ref\/lection
$r_{\gamma}$, one has
\[
r_{\gamma}\delta =\delta -\gamma .
\]
This implies that
\begin{gather*}
(w_{j-1}\delta)(h_{\gamma_j})-1=(r_{j-1}\cdots r_1\delta
)(h_{\gamma_j})-1
=(r_{j-1}\cdots r_2(\delta-\gamma_1))(h_{\gamma_j})-1\\
\phantom{(w_{j-1}\delta)(h_{\gamma_j})-1}{} =(\delta-\gamma_{j-1}-r_{j-1}\gamma_{j-2}-\cdots -r_{j-1}\cdots
r_2\gamma_1)(h_{\gamma_j})-1\\
\phantom{(w_{j-1}\delta)(h_{\gamma_j})-1}{}
=-\sum_{i=1}^{j-1}(r_{j-1}\cdots r_{i+1}\gamma_i)(h_{\gamma_j})=-\sum_{
i=1}^{j-1}\big(w_{j-1}w_{i-1}^{-1}\gamma_i\big )(h_{\gamma_j}).
\end{gather*} This
completes the proof of (b).
\end{proof}

Fix an integral functional $\lambda\in \mathfrak h^{*}$ which is
antidominant. Denote the (algebraic) lowest weight module
corresponding to $\lambda$ by $L(\lambda )$, and a lowest weight
vector by $\sigma_{\lambda}$.  Let $\Phi$ denote the simple roots
$\gamma$ for which $\lambda (h_{\gamma})=0$, where $h_{\gamma}$ is
the coroot, $P=P( \Phi )$ the corresponding parabolic subgroup.
The Borel--Weil theorem in this context realizes $L(\lambda )$ as
the space of strongly regular functions on $G$ satisfying
\[
f(gp)=f(g)\lambda (p)^{-1}
\]
for all $g\in G$ and $p\in P$, where we have implicitly identif\/ied
$\lambda$ with the character of $P$ given by
\[
\lambda (u_1w\exp(x)u_2)=\exp\lambda (x) \] for $x\in \mathfrak
h,\;u_1,u_2\in N^{+},\;w\in W(\Phi )$.  Thus we can view
$L(\lambda )$ as a space of sections of the line bundle
\[ \mathcal L_{\lambda}=G\times_{\lambda}\mathbb
C\rightarrow G/P. \] If $\mathfrak g$ is of f\/inite type, then
$L(\lambda )=H^0(\mathcal L_{\lambda} )$; if $\mathfrak g$ is
af\/f\/ine (and untwisted), then $L(\lambda )$ consists of the
holomorphic sections of f\/inite energy, as in \cite{PS}.

Normalize $\sigma_{\lambda}$ by $\sigma_{\lambda}(1)=1$.

\begin{proposition}\label{dougA26}
 Let $w\in W/W(\Phi )$, and let
$w=r_n\cdots r_1$ be a representative of minimal length~$n$.  Let
$w_j=r_j\cdots  r_1$.  The positive roots mapped to negative roots by
$w$ are given by
\[
\tau_j=w_{j-1}^{-1}(\gamma_j),\qquad 1\leq j\leq n;
\]
let $\lambda_j=-\lambda (h_{\tau_j})$, where $h_{\tau}$ is the
coroot corresponding to $\tau$.  Then
\[
\sigma_{\lambda}^w(r_ni_{\gamma_n}(g_n)\cdots r_1i_{\gamma_1}(g_
1))=
\sigma_{\lambda}\big(w_{n-1}^{-1}i_{\gamma_n}(g_n)w_{n-1}\cdots w_1^{-1}
i_{\gamma_2}(g_2)w_1i_{\gamma_1}(g_1)\big)=\prod_1^na_j^{\lambda_j},
\]
where $g=\begin{pmatrix} a&b\\
c&d\end{pmatrix}\in SL_2$.
\end{proposition}

\begin{proof}[Proof of Proposition~\ref{dougA26}]
The claim about the $\tau_j$ follows from Lemma \ref{dougA8}.
None of these roots lie in $\Delta^{+}(\Phi )$,
by the same argument as follows (\ref{dougA8}).  Thus each $\lambda_j>0$.  It
follows that $\Pi a_j^{\lambda_j}$ is nonzero precisely on the set
$SL_2'\times\cdots\times SL_2'$.

Now $\sigma_{\lambda}^w$, viewed as a section of $\mathcal
L_{\lambda} \rightarrow G/P$, is nonzero precisely on the
$w$-translate of the largest stratum,
\[
w\Sigma_1=wP^{-}P=(P^{-})^wwP.
\]
We claim the intersection of this with $\bar {C}_w$ is $C_w$.  In
one direction \[ C_w=\left(N^{+}\cap (P^{-})^w\right)wP\subseteq
(P^{-})^wwP
\]
by (\ref{dougA10}).  On the other hand $(N^{+}\cap (P^{-})^w)$ is
a closed f\/inite dimensional subgroup of $(P^{-})^w$.  Since
$(P^{-})^w$ is topologically equivalent to $w\Sigma_1$, the limit
points of $C_w$ must be in the complement of $w\Sigma_1$.  This
establishes the other direction.

It now follows from Proposition \ref{dougA4} that
$\sigma_{\lambda}^w$ is also nonzero precisely on
$SL_2'\times\cdots\times SL_2'$, viewed as a function of
$(g_n,\dots ,g_1)$.

Write
\begin{gather}
\sigma_{\lambda}\left(w_{n-1}^{-1}i_{\gamma_n}(g_n)w_{n-1}w_{n-
2}^{-1}i_{\gamma_{n-1}}(g_{n-1})w_{n-2}\cdots
i_{\gamma_1}(g_1)\right )\nonumber\\
\label{dougA29}
\qquad{}=\sigma_{\lambda}\left(i_{\tau_n}(g_n)i_{\tau_{n-1}}(g_{n-1})\cdots
i_{\tau_1}(g_1)\right),
\end{gather}
where $i_{\tau_i}(\cdot )=w_{i-1}^{-1}i_{\gamma_i}(\cdot
)w_{i-1}$.  Because \[ w_{i-1}^{-1}(\gamma_i)>0,
\]
$i_{\tau_j}:SL_2\to G$ is a homomorphism onto the root subgroup
corresponding to $\tau_j$ which is compa\-tible with the canonical
triangular decompositions.

For $g=\left(\begin{matrix} a&b\\
c&d\end{matrix} \right)\in SL_2'$, write $g=LDU$, where
\[
L=\left(\begin{matrix} 1&0\\
ca^{-1}&1\end{matrix} \right),\qquad D=\left(\begin{matrix} a&0\\
0&a^{-1}\end{matrix} \right),\qquad U=\left(\begin{matrix} 1&a^{-1}b\\
0&1\end{matrix} \right).
\]
Then for $(g_j)\in SL_2'\times\cdots\times SL_2'$, (\ref{dougA29})
equals
\begin{gather*}
\sigma_{\lambda}(i_{\tau_n}(L_nD_nU_n)\cdots i_{\tau_1}(L_1D_1U_
1))\\
\qquad{}=\sigma_{\lambda}(i_{\tau_n}(L_nU_n')i_{\tau_{n-1}}(L_{n-1}'U_{
n-1}')\cdots i_{\tau_1}(L_1'U_1')i_{\tau_n}(D_n)\cdots i_{\tau_1}
(D_1))
\\
\qquad{} =\sigma_{\lambda}(i_{\tau_n}(L_nU_n')\cdots i_{\tau_1}(L_1'U_1'
))\Pi a_j^{\lambda_j},
\end{gather*}
where each $L_j'$ ($U_j'$) has the same form as $L_j$ ($U_j$,
respectively).  This follows from the fact that~$H$ normalizes
each $\exp(\mathfrak g_{\pm r})$.

Now each $L_j'U_j'\in SL_2'$, so that $i_{\tau_n}(L_nU_n')\cdots
i_{\tau_1}(L_1'U_1')$ is in $\Sigma_1$.  We now conclude that
\[
\sigma_{\lambda}(i_{\tau_n}(L_nU_n')\cdots i_{\tau_1}(L_1'U_1')
)=1,
\]
by the fundamental theorem of algebra, since this is polynomial
and never vanishes.
\end{proof}

\pdfbookmark[1]{References}{ref}
\LastPageEnding

\end{document}